\documentclass[12pt]{article}

\usepackage{amsfonts}
\usepackage{graphicx}

\def\squarebox#1{\hbox to #1{\hfill\vbox to #1{\vfill}}}
\newcommand{\qed}{\hspace*{\fill}
\vbox{\hrule\hbox{\vrule\squarebox{.667em}\vrule}\hrule}\smallskip}
\newtheorem{teorema}{Theorem}[section]
\newtheorem{lema}[teorema]{Lemma}

\newtheorem{proposicao}[teorema]{Proposition}

\newenvironment{prova}{\noindent {\bf Proof:}}{\hfill $\qed $ \newline}
\newenvironment{obs}{\noindent {\bf Remark:}}{\vspace{0.5cm}}
\newenvironment{exemplo}{\noindent {\bf Example:}}{}

\newcommand{\R}{{\mathbb R}}
\newcommand{\N}{{\mathbb N}}
\newcommand{\Z}{{\mathbb Z}}
\newcommand{\C}{\mathbb{C}}
\newcommand{\F}{\mathbb{F}}

\renewcommand{\t}{{\mathbb T}}

\newcommand{\ad}{{\rm ad}}
\newcommand{\Ad}{{\rm Ad}}

\newcommand{\Int}{{\rm Int}}
\newcommand{\Gl}{{\rm Gl}}

\newcommand{\cl}{{\rm cl} \,}
\renewcommand{\int}{{\rm int} \,}

\newcommand{\g}{\mathfrak{g}}
\renewcommand{\k}{\mathfrak{k}}
\newcommand{\s}{\mathfrak{s}}
\renewcommand{\a}{\mathfrak{a}}

\newcommand{\n}{\mathfrak{n}}

\newcommand{\p}{\mathfrak{p}}

\newcommand{\T}{\Theta}

\newcommand{\ov}[1]{{\overline{#1}}}

\newcommand{\wh}[1]{{\widehat{#1}}}
\newcommand{\wt}[1]{{\widetilde{#1}}}

\newcommand{\til}[1]{{\widetilde{#1}}}

\begin{document}

\title{Jordan decomposition and dynamics \\ on flag manifolds}
\author{Mauro Patr\~{a}o\footnote{Supported by FAPEDF grant
no.??/????}\,
and Lucas Seco\footnote{Supported by FAPESP grant no.??/????}\\
with an appendix by\\ Thiago Ferraiol\footnote{Supported by FAPESP
grant no.??/????}} \maketitle

\begin{abstract}
Let $\g$ be a semisimple Lie algebra and $G = \Int(\g)$. In this
article, we relate the Jordan decomposition of $X \in \g$ (or $g
\in G$) with the dynamics induced on generalized flag manifolds by
the right invariant continuous-time flow generated by $X$ (or the
discrete-time flow generated by $g$). We characterize the
recurrent set and the finest Morse decomposition (including its
stable sets) of these flows and show that its entropy always
vanishes. We characterize the structurally stable ones and compute
the Conley index of the attractor Morse component. When the
nilpotent part of $X$ is trivial, we compute the Conley indexes of
all Morse components. Finally, we consider the dynamical aspects
of linear differential equations with periodic coefficients in
$\g$, which can be regarded as an extension of the dynamics
generated by an element $X \in \g$. In this context, we generalize
Floquet theory and extend the previous results to this case.
\end{abstract}



\section{Introduction}

Let $G$ be a linear group acting differentially on a manifold $F$
and $\g$ its Lie algebra. We denote by $g^t$, $t \in \t = \R$ or
$\Z$, the right invariant continuous-time flow generated by $X \in
\g$ or the discrete-time flow generated by $g \in G$.  More
precisely, when $\t = \R$, we have that $g^t = \exp(tX)$ and, when
$\t = \Z$, we have that $g^t$ is the $t$-iterate of $g$. When
$t=1$ we just write $g=g^1$. Throughout this paper, we call $g^t$
a linear flow. It induces a differentiable flow on $F$ given by
$(t, x) \mapsto g^t x$, where $x \in F$ and $t \in \t = \Z$ or
$\R$. We call these flows linearly induced flows.

Take $G = \Int(\g)$, where $\g$ is a semisimple Lie algebra. In
\cite{dkv} it is considered a continuous-time flow generated by real
semisimple element $H \in \g$ acting on the flag manifolds of $\g$:
they show it is a Morse-Bott gradient flow, describe its fixed point
set and their stable sets.   In \cite{hermann} it is analized a
continuous-time flow generated by an element which is the sum of two
commuting elements of $\g$, one of which induces a gradient vector
field and the other generates a one-parameter group of isometries.
In the context of $\g = {\rm sl}(n,\R)$, the articles \cite{ammar,
batterson} study the discrete-time flow generated by an arbitrary
element $g \in {\rm Sl}(n,\R)$: they characterize the structurally
stable ones.

In this article, we study the dynamics of linearly induced flows
$g^t$, for both continuous and discrete times, acting on a
generalized flag manifold $\F_\T$ of $\g$. This context includes,
for example, volume preserving and Hamiltonian linearly induced
flows acting, respectively, on Grassmanian manifolds and on
Grassmanian of isotropic subspaces, such as the Grassmanian of the
Lagrangian subspaces. This dynamics is related with the Jordan
decomposition of the flow $g^t$, which is defined in terms of the
Jordan decomposition of $X$ or $g$ (see Section
\ref{subsecjordan}). We also consider the dynamical aspects of
linear differential equations with periodic coefficients in $\g$,
which can be regarded as an extension of the dynamics generated by
an element $X \in \g$. In what follows we describe the structure
of this article.  We note that we recover, in this setting, the
results of \cite{pss,msm} about flows in flag bundles with chain
recurrent compact Hausdorff base.

In the preliminaries we recall some notions of Conley theory and of
semisimple Lie theory, proving some useful results.

In Section 3, we recall the Jordan decomposition in $G$ and $\g$ and
show that the flow $g^t$ can be written as a product of commuting
flows $g^t = e^t h^t u^t$, where $e^t$, $h^t$, $u^t$ are linear
flows in $G$ which are called, respectively, the elliptic,
hyperbolic and unipotent components of $g^t$. We finish this section
with a result about the good behavior of the Jordan decomposition
under a certain representation of $G$ which is related to a natural
immersion of a flag manifold into a projective space.

Section 4 is made up of various subsections.  In the first one, we
look at the linearly induced flow of $g^t$ on a flag manifold as the
restriction of a linearly induced flow on a projective space. Using
the results of the appendix about dynamics on projective spaces and
the results of \cite{dkv} about the action of a real semisimple
element $H \in \g$ on the flag manifolds, we generalize these
results characterizing the recurrent and chain recurrent sets, the
finest Morse decomposition, including its stable sets, in terms of
the fixed points of the Jordan components.  For example, we get the
following result.

\begin{teorema}
The recurrent and chain recurrent sets of $g^t$ in the flag manifold
$\F_\T$ are given, respectively, by
\[
\mathcal{R}(g^t) = \emph{fix}_\T\left(h^t\right) \cap
\emph{fix}_\T\left(u^t\right)\quad\mbox{and}\quad \mathcal{R}_C(g^t)
= \emph{fix}_\T\left(h^t\right),
\]
where $\emph{fix}_\T(h^t)$ and $\emph{fix}_\T(u^t)$ are the fixed
points of these flows in $\F_\T$.
\end{teorema}
As a byproduct, we show that the entropy of these flows always
vanishes.  In Section 4.2 we define the conformal flows as the ones
whose unipotent part in the Jordan decomposition is trivial, this is
the kind of linear flow considered in \cite{hermann}.  For these
flows, we compute the Conley indexes of all Morse components.  We
note that we can compute the Conley index of the attractor for every
flow $g^t$, with no restrictions. We then introduce the regular
flows, which are a particular case of the conformal flows. We show
that they are dense in $G$ or $\g$, which implies Theorem 8.1 of
\cite{hermann} about the density of continuous-time conformal flows.
Using this and the previous results we obtain the next result which
generalizes results of \cite{ammar, batterson} obtained in the
context of discrete-time flows generated by an arbitrary element $g
\in {\rm Sl}(n,\R)$.

\begin{teorema}\label{teoestab}
The following conditions are equivalent:
\begin{itemize}
\item[(i)] $g^t$ is regular,

\item[(ii)] $g^t$ is Morse-Smale and

\item[(iii)] $g^t$ is structurally stable.
\end{itemize}
\end{teorema}

Finally, we consider the dynamical aspects of linear differential
equations with periodic coefficients in $\g$, which can be
regarded as an extension of the dynamics generated by an element
$X \in \g$. In this context, we generalize Floquet theory and then
extend the previous results to this case.

\section{Preliminaries}\label{secpreliminar}

\subsection{Flows on topological spaces\label{preliminar-conley}}

Let $\phi :\mathbb{T}\times X\to X$ be a continuous flow on a
compact metric space $(X, d)$, with discrete $\mathbb{T}={\mathbb
Z}$ or continuous $\mathbb{T}={\mathbb R}$ time. For a
$\phi\,^t$-invariant set $ \mathcal{M}\subset X$, we define its
stable and unstable sets respectively as
\[
\mathrm{st}(\mathcal{M})=\{x\in E:\omega (x)\subset
\mathcal{M}\},\qquad \mathrm{un}(\mathcal{M})=\{x\in E:\omega
^{*}(x)\subset \mathcal{M}\},
\]
where $\omega(x)$, $\omega^*(x)$ are the limit sets of $x$. We
denote by $\mathcal{R}(\phi\,^t)$ the set of all recurrent points,
that is
\[
\mathcal{R}(\phi\,^t) = \{x \in X : x \in \omega(x)\},
\]
and by ${\rm fix}(\phi\,^t)$ the set of all fixed points, that is
\[
{\rm fix}(\phi\,^t) = \{x \in X : \phi\,^t(x) = x, \mbox{ for all
} t \in \t\}.
\]

A linear flow $\Phi\,^t$ on a vector bundle $V$ is called normally
hyperbolic if $V$ can be written as a Whitney sum of their stable
and unstable set and there exist a norm in $V$ and constants
$\alpha, \beta > 0$ such that $|\Phi\,^t(v)| < {\rm e}^{-\alpha
t}|v|$, when $v$ is in the stable set, and $|\Phi\,^t(v)| < {\rm
e}^{\beta t}|v|$, when $v$ is in the unstable set. We say that a
$\phi\,^t$-invariant set $ \mathcal{M}\subset X$ is normally
hyperbolic if there exists a neighborhood of $\mathcal{M}$ where
the flow is conjugated to a normally hyperbolic linear flow
restricted to some neighborhood of the null section.

We recall here the definitions and results related to the concept
of chain recurrence and chain transitivity introduced in \cite{c}
(see also \cite{p}). Take $x,y\in X$, $\varepsilon > 0$ and $t\in
\Bbb{T}$. A $( \varepsilon,t)$-chain from $x$ to $y$ is a sequence
of points $ \{x=x_{1},\ldots ,x_{n+1}=y\}\subset X$ and a sequence
of times $\{t_{1},\ldots ,t_{n}\}\subset \Bbb{T}$ such that
$t_{i}\geq t$ and $d(\phi\,^{t_{i}}(x_{i}),x_{i+1}) <
\varepsilon$, for all $i=1,\ldots ,n$.

Given a subset $Y\subset X$ we write $\Omega (Y,\varepsilon,t)$
for the set of all $x$ such that there is a
$(\varepsilon,t)$-chain from a point $y\in Y$ to $x$. Also we put
\[
\Omega ^{*}(x,\varepsilon,t)=\{y\in X:x\in \Omega
(y,\varepsilon,t)\}.
\]
If $Y\subset X$, we write
\[
\Omega (Y)=\bigcap \{\Omega (Y,\varepsilon,t):\varepsilon > 0,t\in
\Bbb{T}\}.
\]
Also, for $x\in X$ we write $\Omega (x)=\Omega(\{x\})$ and define
the relation $x\preceq y$ if $y\in \Omega(x)$, which is
transitive, closed and invariant by $\phi\,^t $, i.e., we have
$\phi\,^{t}(x)\preceq \phi\,^{s}(x)$ if $x\preceq y$, for all
$s,t\in \Bbb{T}$. For every $Y\subset X$ the set $ \Omega (Y)$ is
invariant as well.

Define the relation $x\sim y$ if $x\preceq y$ and $y\preceq x$.
Then $x\in X$ is said to be chain recurrent if $x\sim x$. We
denote by $\mathcal{R}_C(\phi\,^t)$ the set of all chain recurrent
points. It is easy to see that the restriction of $\sim $ to
$\mathcal{R}_C(\phi\,^t)$ is an equivalence relation. An
equivalence class of $\sim $ is called a \textit{chain transitive
component} or a \textit{chain component}, for short.

Now we prove two results which will be used further on.

\begin{lema}\label{lemaiso}
Let $e^t$ be a flow of $X$ such that $e^t$ is an isometry for all
$t \in \t$. Then, for each $T \in \t$ and each $x \in X$, there
exists a sequence $n_k \to \infty$ such that $g^{n_k T}x \to x$.
\end{lema}
\begin{prova}
By the compactness of $X$, we have that the sequence $e^{n T}x$
has a convergent subsequence. Thus, given $\varepsilon
> 0$ and $L > 0$, there exist $m,k \in \N$ such that $m - k > L$
and
$$
d(e^{(m-k)T}x,\, x) = d(e^{mT}x,\, e^{kT}x) < \varepsilon.
$$
Hence there exists a sequence $n_k\to \infty$ such that $g^{n_k
T}x \to x$.
\end{prova}

\begin{lema}\label{lemaeu}
Let $e^t$, $u^t$ be commuting flows of $X$, $t \in \t$. Assume
that $e^t$ is an isometry for all $t \in \t$ and that for each $x
\in X$ there exists $y \in X$ such that the omega and alpha limits
of $x$ by $u^t$ are precisely $y$. Then the composition $e^tu^t$
is a chain recurrent flow.
\end{lema}
\begin{prova}
Fix $x \in X$.  Given $\varepsilon > 0$ and $t_0
> 0$ we will construct an $(\varepsilon,t_0)$-chain from $x$ to
$x$. By the assumption on $u$, there exists $y \in X$ and $t_1
> t_0$ such that
$$
u^t (x),\, u^{-t}(x) \in B(y,\varepsilon/2),
$$
for all $t > t_1$. Taking $t > t_1$, it follows that the points
$\{x, u^{-t}(x), x\}$ and times $\{t,t\}$ define an
$(\varepsilon,t_0)$-chain of $u$, since
$$
d(u^t(x),\, u^{-t}(x)) < \varepsilon \qquad\mbox{and}\qquad d(u^t
u^{-t} (x),x)=0 < \varepsilon.
$$
Now, since the isometry $e$ is recurrent (see Lemma
\ref{lemaiso}), there exists $t > t_1$ such that $d(e^{2t}(x), x)
< \varepsilon$.  Thus the points $\{x, e^t u^{-t} (x), x\}$ and
times $\{t,t\}$ define an $(\varepsilon,t_0)$-chain of $eu$. In
fact, using the commutativity of $e$ and $u$ and using that $e$ is
an isometry, we have
$$
d( (eu)^t (x), e^t u^{-t} (x) ) = d(u^t (x),\, u^{-t}(x)) <
\varepsilon,
$$
by the above construction. Finally, using again the commutativity
of $e$ and $u$ we have that
$$
d( (eu)^t e^t u^{-t} (x), x ) = d(e^{2t}(x), x) < \varepsilon,
$$
by the choice of $t$.
\end{prova}

Now we relate Morse decompositions to chain transitivity. First
let us recall that a finite collection of disjoint subsets $
\{{\cal M}_{1},\ldots ,{\cal M}_{n}\}$ defines a Morse
decomposition when
\begin{itemize}
\item[(i)] each ${\cal M}_i$ is compact and $\phi\,^t$-invariant,

\item[(ii)] for all $x \in X$ we have $\omega(x),\, \omega^*(x)
\subset \bigcup_i {\cal M}_i$,

\item[(iii)] if $\omega(x),\, \omega^*(x) \subset {\cal M}_j$ then
$x \in {\cal M}_j$.
\end{itemize}
Each set ${\cal M}_{i}$ of a Morse decomposition is called a Morse
component. If $ \{{\cal M}_{1},\ldots ,{\cal M}_{n}\}$ is a Morse
decomposition of $X$, then it is immediate that $X$ decomposes as
the disjoint union of stable sets $\mathrm{st}(\mathcal{M}_{i})$.

The finest Morse decomposition is a Morse decomposition which is
contained in every other Morse decomposition. The existence of a
finest Morse decomposition of a flow is equivalent to the
finiteness of the number of chain components (see \cite {p},
Theorem 3.15).  In this case, each Morse component is a chain
transitive component and vice-versa. We say that the flow
$\phi\,^t$ is normally hyperbolic if there exists the finest Morse
decomposition and their Morse components are normally hyperbolic.

\subsection{Semi-simple Lie theory\label{preliminar-lie}}

For the theory of semisimple Lie groups and their flag manifolds
we refer to Duistermat-Kolk-Varadarajan \cite{dkv}, Helgason
\cite{helgason} and Warner \cite{w}. To set notation let $\g$ be a
semisimple Lie algebra and $G = \Int(\g) \subset {\rm Gl}(\g)$
acting in $\g$ canonically. We identify throughout the Lie algebra
of $G$ with $\g$, that is, we write $g = \exp(X)$ to mean $g =
{\rm e}^{\ad(X)}$, where $X \in \g$. Thus, for $g \in G$ and $X
\in X$, it follows that $g \exp(X) g^{-1} = \exp( gX )$.  Note
that if $\wt{G}$ is a connected Lie group with Lie algebra $\g$,
then $\Ad(\wt{G}) = \Int(\g) = G$.  It follows that the adjoint
action of $\wt{G}$ in $\g$ is the canonical action of $G$ in $\g$.

Fix a Cartan involution $\theta $ of $\frak{g}$ with Cartan
decomposition $\frak{g}=\frak{k}\oplus \frak{s}$. The form
$B_{\theta }\left( X,Y\right) =-\langle X,\theta Y\rangle $, where
$\langle \cdot ,\cdot \rangle $ is the Cartan-Killing form of
$\frak{g}$, is an inner product.

Fix a maximal abelian subspace $\frak{a}\subset \frak{s}$ and a
Weyl chamber $\frak{a}^{+}\subset \frak{a}$. We let $\Pi $ be the
set of roots of $\frak{a }$, $\Pi ^{+}$ the positive roots
corresponding to $\frak{a}^{+}$, $\Sigma $ the set of simple roots
in $\Pi ^{+}$ and $\Pi^- = - \Pi^+$ the negative roots. The
Iwasawa decomposition of the Lie algebra $\frak{g}$ reads
$\frak{g}= \frak{k}\oplus \frak{a}\oplus \frak{n}^{\pm}$ with
$\frak{n} ^{\pm}=\sum_{\alpha \in \Pi ^{\pm}}\frak{g}_{\alpha }$ where
$\frak{g} _{\alpha }$ is the root space associated to $\alpha $.
As to the global decompositions of the group we write $G=KS$ and
$G=KAN^\pm$ with $K=\exp \frak{k}$, $S=\exp \frak{s}$, $A=\exp
\frak{a}$ and $N^{\pm}=\exp \frak{n} ^{\pm}$.

The Weyl group $W$ associated to $\frak{a}$ is the finite group
generated by the reflections over the root hyperplanes $\alpha =0$
in $\frak{ a}$, $\alpha \in \Pi $. For each $w \in W$ and $\alpha
\in \Pi$ we define $w^*\alpha(H) = \alpha(w^{-1}H)$, for all $H \in
\a$. We have that $w^*\alpha \in \Pi$ and that this is a transitive
action of $W$ on $\Pi$.  The maximal involution $w^-$ of $W$ is the
(only) element of $W$ which is such that $(w^-)^*\Sigma = - \Sigma$.

Given a subset of simple roots $\Theta \subset \Sigma $, let
$$
\frak{ a}_{\Theta }=\{H\in \frak{a}:\alpha (H)=0,\,\alpha \in
\Theta \}
$$
and put $A_\T = \exp(\a_\T)$. Let also
$$
\frak{n}({\Theta })^{\pm }=\sum_{\alpha \in \langle \Theta \rangle
\cap \Pi^\pm }\frak{g}_{\alpha } \qquad\mbox{and}\qquad
\frak{n}_{\Theta }^{\pm }=\sum_{\alpha \in \Pi ^{\pm }-\langle
\Theta \rangle }\frak{g}_{\alpha }
$$
and put $N_{\Theta }^{\pm }=\exp (\frak{n}_{\Theta }^{\pm })$. The
subset $\Theta $ singles out the subgroup $W_{\Theta }$ of the
Weyl group which acts trivially on $\frak{a}_{\Theta }$.

The standard parabolic subalgebra of type $\Theta \subset \Sigma$
with respect to chamber $\frak{a}^+$ is defined by
\[
\frak{p}_{\Theta }=\frak{n}^{-}\left( \Theta \right) \oplus
\frak{m}\oplus \frak{a}\oplus \frak{n}^{+}.
\]
Let $p$ the dimension of $\p_\T$ and denote the grassmanian of
$p$-dimensional subspaces of $\g$ by ${\rm Gr}_p(\g)$. The flag
manifold of type $\Theta $ is the orbit $\Bbb{F}_{\Theta
}=G\frak{p}_{\Theta } \subset {\rm Gr}_p(\g)$, with base point
$b_{\Theta }=\frak{p}_{\Theta }$, which identifies with the
homogeneous space $G/P_{\Theta }$. Since the center of $G$
normalizes $\frak{p}_{\Theta }$, the flag manifold depends only on
the Lie algebra $\frak{g}$ of $G$. The empty set $\Theta
=\emptyset $ gives the maximal flag manifold
$\Bbb{F}=\Bbb{F}_{\emptyset }$ with basepoint $b=b_{\emptyset }$.

For $ H\in \frak{a}$ we denote by $Z_{H}$, $K_H$, $W_{H}$, the
centralizer of $H$ in $G$, $K$, $W$, respectively, i.e, the
elements in those groups which fix $H$. Note that $g \in G$
centralizes $H$ if and only if it commutes with $\ad(H)$. In fact,
this follow from $\ad(gH) = g\ad(H)g^{-1}$ and the injectivity of
$\ad$. When $H\in \mathrm{cl}\frak{ a}^{+}$ we put
\[
\Theta (H)=\{\alpha \in \Sigma :\alpha (H)=0\}.
\]

An element $H\in \mathrm{cl}\frak{a}^{+}$ induces a vector field $
\widetilde{H}$ on a flag manifold ${\mathbb F}_{\Theta }$ with
flow $\exp (tH)$. This is a gradient vector field with respect to
a given Riemannian metric on ${\mathbb F}_{\Theta }$ (see
\cite{dkv}, Section 3). The connected sets of fixed point of this
flow are given by
\[
\mathrm{fix}_{\Theta }(H,w)=Z_{H} wb_{\Theta }=K_{H} wb_{\Theta },
\]
so that they are in bijection with the cosets in $W_{H}\backslash
W/W_{\Theta }$. Each $w$-fixed point connected set has stable
manifold given by
\[
\mathrm{st}_{\Theta }(H,w)=N_{\T(H)}^{-} \mathrm{fix}_{\Theta
}(H,w),
\]
whose union gives the Bruhat decomposition of ${\mathbb F}_{\Theta
}$:
\[
{\mathbb F}_{\Theta }=\coprod_{W_{H}\backslash W/ W_{\Theta
}}\mathrm{st}_{\Theta }(H,w).
\]
The unstable manifold is
\[
\mathrm{un}_{\Theta }(H,w)=N_{\T(H)}^{+} \mathrm{fix}_{\Theta
}(H,w).
\]
We note that both $\mathrm{st}_{\Theta }(H,1)$ and
$\mathrm{un}_{\Theta }(H,w^-)$ are open and dense in $\F_\T$. Since
the centralizer $Z_{H}$ of $H$ leaves $\mathrm{fix}_{\Theta }(H,w)$
invariant and normalizes both $N_{\T(H)}^{-}$ and $N_{\T(H)}^{+}$,
it follows $ \mathrm{st}_{\Theta }(H,w)$ and $\mathrm{un}_{\Theta
}(H,w)$ are $Z_{H}$ -invariant. We note that these fixed points and
(un)stable sets remain the same if $H$ is replaced by some
$H^{\prime }\in \mathrm{cl}\frak{a}^{+}$ such that $\Theta (H') =
\Theta(H)$.

We note that, since the spectrum of $\ad(H)$ in $\g_\alpha$ is
$\alpha(H)$, it follows that the spectrum of $h = \exp(H)$ in
$\g_\alpha$ is $e^{\alpha(H)}$.

We conclude with a useful lemma about the decomposition semisimple
elements.  We say that $X \in \g$ is semisimple if $\ad(X)$ is
diagonalizable over $\C$ and that $ g \in G$ is semisimple if $g$
is diagonalizable over $\C$.

\begin{lema}\label{lemasemisimples}
We have that
\begin{itemize}
\item[(i)] If $X \in \g$ is semisimple, then there exists an
Cartan decomposition $\g = \k \oplus \s$ such that $X = E + H$
where $H \in \s$ and $E \in \k_H$.

\item[(ii)]  If $g \in G$ is semisimple, then there exists an
Cartan decomposition $G = KS$ such that $g =eh$, where $h =
\exp(H)$, $H \in \s$ and $e \in K_H$.
\end{itemize}
\end{lema}
\begin{prova}
For item (i), since $X$ is semisimple, there exists a Cartan
subalgebra $\frak{j}$ such that $X \in \frak{j}$ (see the proof of
Proposition 1.3.5.4, p.105 of \cite{w}).  By Proposition 1.3.1.1,
p.89 of \cite{w}, there exists a Cartan involution $\theta$ such
that $\frak{j}$ is $\theta$-invariant. Thus we have that
\[
\frak{j}=(\frak{j}\cap \k)\oplus (\frak{j}\cap\frak{s}).
\]
Writing $X=E + H$, with $E \in \frak{j} \cap \k$ and $H \in
\frak{j}\cap\frak{s}$, we have that $E$ and $H$ commute, since
$\frak{j}$ is abelian.

For item (ii), since $g$ is semisimple, there exists a Cartan
subgroup $J$ such that $g \in J$ (since the centralizer of $g$ in
$\g$ contains a Cartan subalgebra, see the proof of Proposition
1.4.3.2, p.120 of \cite{w}).  Denote by $\frak{j}$ the associated
Cartan subalgebra. By Proposition 1.3.1.1, p.89 of \cite{w}, there
exists a Cartan involution $\theta$ such that $\frak{j}$ is
$\theta$-invariant. Thus, by Proposition 1.4.1.2, p.109 of
\cite{w}, we have that
\[
J=(J\cap K)(\exp(\frak{j}\cap\frak{s})).
\]
Writing $g=eh$, with $e \in J\cap K$ and $h = \exp(H)$, where $H
\in \frak{j}\cap\frak{s}$. Since $J$ centralizes $\frak{j}$, it
follows that $e$ and $\ad(H)$ commute, showing that $e \in K_H$.
\end{prova}

\section{Jordan decomposition}\label{subsecjordan}

In this section we recall the additive and the multiplicative
Jordan decompositions. Let $V$ be a finite dimensional vector
space.

If $X \in {\rm gl}(V)$, then we can write $X = E + H + N$, where
$E \in {\rm gl}(V)$ is semisimple with imaginary eigenvalues, $H
\in {\rm gl}(V)$ is diagonalizable in $V$ with real eigenvalues
and $N \in {\rm gl}(V)$ is nilpotent. The linear maps $E$, $H$ and
$N$ commute, are unique and called, respectively, the
\emph{elliptic}, the \emph{hyperbolic}, and the \emph{nilpotent}
components of the additive Jordan decomposition of $X$ (see
Section ? of \cite{humphreys}).

If $g \in {\rm Gl}(V)$, then we can write $g = ehu$, where $e \in
{\rm Gl}(V)$ is an isometry relative to some appropriate inner
product, $h \in {\rm Gl}(V)$ is diagonalizable in $V$ with
positive eigenvalues and $u \in {\rm Gl}(V)$ is the exponential of
a nilpotent linear map. The linear maps $e$, $h$ and $u$ commute,
are unique and called, respectively, the \emph{elliptic}, the
\emph{hyperbolic} and the \emph{unipotent} components of the
multiplicative Jordan decomposition of $g$  (see Lemma IX.7.1
p.430 of \cite{helgason}). We denote by $\log h$ the matrix given
by the logarithm of the diagonal elements of $h$ in the Jordan
basis. Writing $g$, $e$, $h$, $u$ in the Jordan basis, we see that
they commute with $\log h$.

Take $\g$ a semisimple Lie algebra. We say that $X = E + H + N$,
where $E, H, N \in \g$, is the Jordan decomposition of $X$ in $\g$
if $\ad(X) = \ad(E) + \ad(H) + \ad(N)$ is the additive Jordan
decomposition of $\ad(X)$ in ${\rm gl}(\g)$. In this case, $E$,
$H$ and $N$ commute, are unique and called, respectively, the
\emph{elliptic}, the \emph{hyperbolic}, and the \emph{nilpotent}
components of $X$.

We note that the conjugate of a Jordan decomposition is the Jordan
decomposition of the conjugate. Now we prove the following useful
result.

\begin{lema}\label{lemajordan}
Let $G = \int(\g)$, where $\g$ is a semisimple Lie algebra. Then
we have that
\begin{itemize}
\item[(i)] For each $X \in \g$, there exists the Jordan
decomposition $X = E + H + N$. Furthermore, there exists an
Iwasawa decomposition $\g = \k \oplus \a \oplus \n^+$ such that $E
\in \k_H$ and $H \in \cl \a^+$.

\item[(ii)] For each $g \in G$, its multiplicative Jordan
components $e, h, u$ lie in $G$. Moreover, there exist a unique $H
\in \g$ such that $\log h = \ad(H)$ and an Iwasawa decomposition
$G = KAN$ such that $e \in K_H$ and $H \in \cl \a^+$.
\end{itemize}
\end{lema}
\begin{prova}
For item (i), by Proposition 1.3.5.1, p.104 of \cite{w}, there
exists a unique decomposition $X = S + N$, where $S, N \in \g$
commute, $S$ is semisimple and $\ad(N)$ is nilpotent. By Lemma
\ref{lemasemisimples}, there exists an Cartan decomposition $\g =
\k \oplus \s$ such that $S = E + H$, where $H \in \s$ and $E \in
\k_H$. This is the additive Jordan decomposition of $S$ in $\g$,
since $\ad(E)$ is $B_{\theta}$-anti-symmetric and $\ad(H)$ is
$B_{\theta}$-symmetric. It remains to show that $\ad(E)$ and
$\ad(H)$ commute with $\ad(N)$. We first note that $u = I +
\ad(N)$ is invertible and that $Y \in {\rm gl}(\g)$ commutes with
$\ad(N)$ if and only if $Y$ commutes with $u$. In fact, we have
that
\[
Y + Y\ad(N) = Yu = uY = Y + \ad(N)Y.
\]
It follows that $u$ commutes with $\ad(X)$. In order to show that
$\ad(E)$ and $\ad(H)$ commute with $u$, we write
\[
\ad(E) + \ad(H) + \ad(N) = u\ad(X)u^{-1} = u\ad(E)u^{-1} +
u\ad(H)u^{-1} + \ad(N).
\]
By the uniqueness of the additive Jordan decomposition in ${\rm
gl}(\g)$, we have that $\ad(E) = u\ad(E)u^{-1}$ and $\ad(H) =
u\ad(H)u^{-1}$. Since $H \in \s$, we can choose an Iwasawa
decomposition $\g = \k \oplus \a \oplus \n^+$ such that $E \in
\k_H$ and $H \in \cl \a^+$.

For item (ii), by Proposition 1.4.3.3, p.120 of \cite{w}, there
exists a unique decomposition $g = \wh{s} \wh{u}$, where $\wh{s} ,
\wh{u} \in G$ commute, $\wh{s} $ is semisimple and $\wh{u}$ is the
exponential of a nilpotent linear map. By Lemma
\ref{lemasemisimples}, there exists an Cartan decomposition $G = KS$
such that $\wh{s} = \widehat{e}\widehat{h}$, where $\wh{h} =
\exp(\wh{H})$, $\wh{H} \in \s$ and $\wh{e} \in K_{\wh{H}}$. This is
the multiplicative Jordan decomposition of $g$, since $\wh{e}$ is a
$B_{\theta}$-isometry and $\wh{h}$ is $B_{\theta}$-positive. In
order to show that $\wh{e}$ and $\wh{h}$ commute with $\wh{u}$, one
can proceed as in Lemma IX.7.1 p.431 of \cite{helgason}. By the
uniqueness of the multiplicative Jordan decomposition in $\Gl(\g)$,
it follows that $e = \wh{e}$, $h = \wh{h}$ and $u = \wh{u}$, showing
that the multiplicative Jordan components of $g$ lie in $G$. By the
proof of Lemma IX.7.3 item (i) p.431 of \cite{helgason}, we have
that $\log h$ lies in the Lie algebra of $G$ and thus there exists a
unique $H \in \g$ such that $\log h = \ad(H)$, since $\ad$ is
injective. Since both $\ad(H)$ and $\ad(\wh{H})$ commute with $e =
\wh{e}$, it follows that $\ad(H)$ and $\ad(\widehat{H})$ can be
diagonalized in the same basis. Since ${\rm e}^{\ad(H)} = h = {\rm
e}^{\ad(\wh{H})}$ and using the injectivity of $\ad$, it follows
that $H = \widehat{H} \in \s$. Thus we can choose an Iwasawa
decomposition $G = KAN$ such that $e \in K_H$ and $H \in \cl \a^+$.
\end{prova}

Let $G$ be a linear group.  Now we define the Jordan decomposition
of a linear flow $g^t$ in $G$, $t \in \t$. If $\t = \R$ then $g^t
= {\rm e}^{tX}$, $X \in {\rm gl}(V)$ and we can use the additive
Jordan decomposition $X = E + H + N$ to write $g^t = e^t h^t u^t$,
where $e^t = {\rm e}^{tE}$, $h^t = {\rm e}^{tH}$ and $u^t = {\rm
e}^{tN}$.  If $\t = \Z$ we can use the multiplicative Jordan
decomposition to write $g^t = e^t h^t u^t$ for each $t \in \t$. It
follows that in both cases the linear flows $g^t$, $e^t$, $h^t$,
$u^t$ commute.

Now take $G=\Int(\g)$, where we identify the Lie algebra of $G$
with $\g$ (see Section \ref{preliminar-lie}).  Let $g^t \in G$,
for all $t \in \t$.  By Lemma \ref{lemajordan}, each Jordan
component $e^t, h^t, u^t$ of $g^t$ also lies in $G$, for all $t
\in \t$.  If $\t = \Z$ then this is immediate. When $\t = \R$,
then $g^t = \exp(tX)$, where $X \in \g$. Thus we can use the
Jordan decomposition $X = E + H + N$ to write $g^t = e^t h^t u^t$,
where $e^t = \exp(tE)$, $h^t = \exp(tH)$ and $u^t = \exp(tN)$. In
both continuous and discrete time cases, we also have that each
Jordan components of the flow $g^t$ lie in $Z_H$, where $H$ is
given by Lemma \ref{lemajordan}, when $\t = \Z$.

Let $\rho : G \to {\rm Gl}(V)$ be a finite dimensional
representation, where ${\rm d}_1\rho : \g \to {\rm gl}(V)$ is its
infinitesimal representation. When $t \in \Z$, it is immediate
that $\rho(g^t) = \rho(g)^t$. When $t \in \R$, we have that $g^t =
\exp(tX)$, for $X \in \g$. Denoting $\rho(g)^t = \exp(t{\rm
d}_1\rho X)$, it follows also that $\rho(g^t) = \rho(g)^t$.

Now we consider the behavior of the Jordan decomposition with
respect to the canonical representation of the general linear
group ${\rm Gl}(L)$ in ${\rm Gl}(\bigwedge^p L)$, where $L$ be a
finite dimensional vector space, given by
\[
\rho(g) v_1 \wedge \cdots  \wedge v_p = gv_1 \wedge \cdots  \wedge
gv_p.
\]

\begin{lema}\label{lemarepresent}
For the the canonical representation of the general linear group
${\rm Gl}(L)$ in ${\rm Gl}(\bigwedge^p L)$, we have that
\begin{itemize}
\item[i)]  Take $X \in {\rm gl}(L)$.  If $X$ is elliptic (resp.\
hyperbolic, nilpotent), then ${\rm d}_1 \rho X$ is elliptic
(resp.\ hyperbolic, nilpotent). In particular, if $X = E + H + N$
is the additive Jordan decomposition of $X$, then ${\rm d}_1 \rho
X = {\rm d}_1 \rho E + {\rm d}_1 \rho  H + {\rm d}_1 \rho N$ is
the additive Jordan decomposition of ${\rm d}_1 \rho X$.

\item[ii)] Take $g \in {\rm Gl}(L)$.  If $g$ is elliptic (resp.\
hyperbolic, unipotent), then $\rho(g)$ is elliptic (resp.\
hyperbolic, unipotent). In particular, if $g = ehu$ is the Jordan
decomposition of $g$, then $\rho(g) = \rho(e) \rho(h) \rho(u)$ is
the Jordan decomposition of $\rho(g)$.

\item[iii)] If $e^t h^t u^t$ is the Jordan decomposition of $g^t$,
then $\rho(e)^t \rho(h)^t \rho(u)^t$ is the Jordan decomposition
of $\rho(g)^t$.
\end{itemize}
\end{lema}
\begin{prova}
First observe that the complexification of wedge product of $L$ is
equal to the wedge product of the complexification of $L$, that is
$(\bigwedge^p L)_{\C} = \bigwedge^p L_{\C}$.  In fact, it is
immediate that $\bigwedge^p L_{\C} \subset (\bigwedge^p L)_{\C}$
and that both have the same dimension. Note also that
$$
({\rm d}_1 \rho X) v_1 \wedge \cdots  \wedge v_p = \sum_{i = 1}^p
v_1 \wedge \cdots  \wedge X v_i \wedge \cdots \wedge v_p.
$$
We claim that
$$
({\rm d}_1 \rho X)^m v_1 \wedge \cdots  \wedge v_p = \sum_i w_{i1}
\wedge \cdots \wedge w_{ip},
$$
where $w_{ij} \in N^{m_{ij}}(L)$ such that $\sum_{j=1}^p m_{ij} =
m$. For $m=0$ this is immediate.  By induction on $m$
$$
({\rm d}_1 \rho X)^{m+1} v_1 \wedge \cdots  \wedge v_p = ({\rm
d}_1 \rho X) \sum_i w_{i1} \wedge \cdots \wedge w_{ip} =
$$
$$
= \sum_i \sum_{j = 1}^p  w_{i1} \wedge \cdots  \wedge X w_{ij}
\wedge \cdots  \wedge w_{ip}
$$
Since $w_{ij} \in N^{m_{ij}}(L)$ it follows that $Xw_{ij} \in
N^{m_{ij}+1}(L)$.

For item (i), taking $X$ nilpotent then there exists $l$ such that
$X^l = 0$.  From the above claim, it follows that $({\rm d}_1 \rho
X)^{pl} = 0$.  In fact,
$$
({\rm d}_1 \rho X)^{pl} v_1 \wedge \cdots  \wedge v_p = \sum_i
w_{i1} \wedge \cdots \wedge w_{ip},
$$
where $w_{ij} \in N^{m_{ij}}(L)$ such that $\sum_{j=1}^p m_{ij} =
pl$.  Thus, for each $i$ there exists $j$ such that $m_{ij} \geq l$.
Therefore $w_{ij} = 0$, which implies that $w_{i1} \wedge \cdots
\wedge w_{ip} = 0$, for all $i$.  Now taking $X$ elliptic, there
exists a $\C$-basis $v_1, \ldots v_n$ of $L_\C$ such that $X v_k =
z_k v_k$, where $z_k$ is purely imaginary.  Then
$$
\{ v_{i_1} \wedge \cdots  \wedge v_{i_p} :\mbox{ where } 1 \leq
i_1 < \ldots < i_p \leq n\}
$$
is a $\C$-basis of $\bigwedge^p L_{\C}$ such that
$$
({\rm d}_1 \rho X) v_{i_1} \wedge \cdots  \wedge v_{i_p} =
(z_{i_1} + \cdots + z_{i_p}) v_{i_1} \wedge \cdots \wedge v_{i_p}.
$$
This implies that $ {\rm d}_1 \rho X $ is elliptic, since $z_{i_1}
+ \cdots + z_{i_p}$ is purely imaginary.  The hyperbolic case is
analogous.

For item (ii), taking $g$ unipotent, then $g = e^N$ with $N$
nilpotent so that $\rho(g) = e^{{\rm d}_1 \rho(N) }$ is unipotent,
by using item (i). Now taking $g$ is elliptic, there exists a
$\C$-basis $v_1, \ldots v_n$ of $L_\C$ such that $g v_k = z_k v_k$,
where $z_k \in \C$ with $|z_k| = 1$.  Then
$$
\{ v_{i_1} \wedge \cdots  \wedge v_{i_p} :\mbox{ where } 1 \leq
i_1 < \ldots < i_p \leq n\}
$$
is a $\C$-basis of $\bigwedge^p L_{\C}$ such that
$$
\rho(g) v_{i_1} \wedge \cdots  \wedge v_{i_p} = (z_{i_1} \cdots
z_{i_p}) v_{i_1} \wedge \cdots  \wedge v_{i_p}.
$$
This implies that $ \rho(g) $ is elliptic, since $|z_{i_1} \cdots
z_{i_p}| = 1$. The hyperbolic case is analogous.

Item (iii) follows immediately from the previous items.
\end{prova}

We recall the well known Pl\"{u}cker embedding, which is given by
\[
i: {\rm Gr}_p(L) \to \mathbb{P}{\textstyle (\bigwedge^p L)}, \quad
P \mapsto [v_1 \wedge \cdots \wedge v_p],
\]
where $\{v_1,\ldots,v_p\}$ is a basis of $P$. This embedding has
the following equivariance property
\[
i(gP) = \rho(g)i(P)
\]
where $\rho$ is the canonical representation presented in Lemma
\ref{lemarepresent}.  If $g^t$ is a linearly induced flow it
follows that
\[
i(g^tP) = \rho(g)^ti(P).
\]

\section{Dynamics in flag manifolds}

In this section, we relate the Jordan decomposition of $g^t$ in $G =
\Int(\g)$ to the dynamics of the induced linear flow $g^t$ on the
flag manifolds of $\g$, where $\g$ is a semisimple Lie algebra. The
main results of the section deals with the characterization of the
recurrent set and the finest Morse decomposition in terms of the
fixed points of the Jordan components.

Recall that, as seen in Section \ref{subsecjordan}, when $g^t \in
G$ then each multiplicative Jordan component $e^t, h^t, u^t$ of
$g^t$ lies in $Z_H$, where $H \in \g$ is such that $\log h =
\ad(H)$. Furthermore, there exists a Weyl chamber $\a^+$ such that
$H \in \cl \a^+$.

Let $\T \subset \Sigma$.  It follows that $e^t$, $h^t$ and $u^t$
induce flows in the flag manifold $\F_\T$. If $p = \dim(\p_\T)$,
we know that $\F_\T \subset {\rm Gr}_p(\g)$, so we can restrict
the Pl\"{u}cker embedding (see Section \ref{subsecjordan}) to $\F_\T$
and get an embedding $i: \F_\T \to \mathbb{P}V$, where $V =
\bigwedge^p \g$. Since $\F_\T$ is $G$-invariant, we have the
following equivariance property
\[
i(gx) = \rho(g)i(x), \quad x \in \F_\T,
\]
where $\rho: G \to {\rm Gl}(V)$ is the restriction to $G$ of the
canonical representation presented in Lemma \ref{lemarepresent}.

\subsection{Recurrence, chain recurrence and entropy}

The next proposition shows that the fixed points of the hyperbolic
part of $g^t$ are Morse components for the flow $g^t$ in the flag
manifold.  This result is proved by using Proposition
\ref{propmorse}, which is in fact a particular case when the flag
manifold is the projective space.

\begin{proposicao}\label{propmorseflag}
Let $g^t$ be a flow on $\F_\T$. The set
$$
\{{\rm fix}_\T(H,w): w \in W_H\backslash W/W_\T \}
$$
is a Morse decomposition for $g^t$.
Furthermore, the stable and unstable sets of ${\rm fix}_\T(H,w)$
are given by
$$
{\rm st} ({\rm fix}_\T(H,w)) = {\rm st}_\T(H,w) \qquad \emph{and}
\qquad {\rm un} ({\rm fix}_\T(H,w)) = {\rm un}_\T(H,w),
$$
\end{proposicao}
\begin{prova}
Since ${\rm fix}_\T(H,w) = Z_H w \p_\T$ and since $g^t \in Z_H$,
it follows that ${\rm fix}_\T(H,w)$ is $g^t$-invariant.

Now we show that ${\rm st}_\T(H,w)$ is the stable set of ${\rm
fix}_\T(H,w)$. By the Bruhat decomposition of $\F_\T$, it is
enough to show that ${\rm st}_\T(H,w)$ is contained in the stable
set of ${\rm fix}_\T(H,w)$. Let $x \in {\rm st}_\T(H,w)$, then $x
= \exp(Y)lwb_\T$, where $Y \in \n^-_{\T(H)}$ and $l \in Z_H$. Then
$$
g^t x = g^t\exp(Y)g^{-t} g^t lw b_\T = \exp(g^t Y)g^t lw b_\T,
$$
where $g^t lw b_\T \in {\rm fix}_\T(H,w)$, since $g^t l \in Z_H$.
Now we show that $g^t Y \to 0$.  This follows by Lemma
\ref{lemalucas1}, since the spectral radius of the restriction of
$g$ to $\n^-_{\T(H)}$ is smaller than 1. In fact, by the Jordan
decomposition, $r(g)$ is given by the greatest eigenvalue of its
hyperbolic component, which is given by the restriction of $h$ to
$\n^-_{\T(H)}$. These eigenvalues are given by $e^{-\alpha(H)}$,
where $\alpha \in \Pi^+$ with $\alpha(H) > 0$, so that $r(g) < 1$.
Now if $g^{t_j}x \to y$ then $g^{t_j} lw b_\T \to y$, so that $y$
lies in the closed subset ${\rm fix}_\T(H,w)$.

For the unstable set we proceed analogously. It follows that
\[
{\rm fix}_\T(h^t) = \bigcup_{w \in W} {\rm fix}_\T(H,w)
\]
contains all the alpha and omega limit sets. In order to show that
the set $\{{\rm fix}_\T(H,w): w \in W_H\backslash W/W_\T \}$ is a
Morse decomposition for $g^t$ it is enough to prove that if
$\omega(x)$, $\omega^*(x) \subset {\rm fix}_\T(H,w)$, then $x \in
{\rm fix}_\T(H,w)$. First recall that
$$
i( {\rm fix}_\T(h^t) ) = {\rm fix}(\rho(h)^t) \cap i(\F_\T),
$$
$$
i( \omega(x) ) = \omega(i(x))\qquad\mbox{and}\qquad i( \omega^*(x)
) = \omega^*(i(x)),
$$
where $i$ is the Pl\"{u}cker embedding. By hypothesis $\omega(i(x))$,
$\omega^*(i(x))$ are contained in the connected set $i({\rm
fix}_\T(H,w))$ of ${\rm fix}(\rho(h)^t)$, so they lie in the same
connected component of ${\rm fix}(\rho(h)^t)$ which is given by an
eigenspace of $\rho(h)$. Using Lemmas \ref{lemarepresent} and
\ref{lemthiago}, it follows that $i(x) \subset {\rm
fix}(\rho(h)^t)$, which shows that $x \in {\rm fix}_\T(h^t)$. Then
there exists $s \in W$ such that $x \in {\rm fix}_\T(H,s)$. By the
invariance of ${\rm fix}_\T(H,s)$, we get that
\[
\omega(x) \subset {\rm fix}_\T(H,s) \cap {\rm fix}_\T(H,w),
\]
showing that $x \in {\rm fix}_\T(H,w)$. The proof for the unstable
set is completely analogous.
\end{prova}

We note that ${\rm st} ({\rm fix}_\T(H,1))$ and ${\rm un} ({\rm
fix}_\T(H,w^-))$ are open and dense (see Section
\ref{preliminar-lie}) so that ${\rm fix}_\T(H,1)$ and ${\rm
fix}_\T(H,w^-)$ are, respectively, the only attractor and repeller
which are thus denoted by ${\rm fix}^+_\T(H)$ and ${\rm
fix}^-_\T(H)$.  Using the previous result, we obtain the desired
characterization of the finest Morse decomposition.

\begin{teorema}\label{teomorseflag}
Let $g^t$ be a flow on $\F_\T$ and $g^t = e^th^tu^t$ its Jordan
decomposition. Each ${\rm fix}_\T(H,w)$ is chain transitive, so
that $\{{\rm fix}_\T(H,w): w \in W_H\backslash W/W_\T \}$ is the
finest Morse decomposition. In particular, the chain recurrent set
of $g^t$ in $\F_\T$ is given by
\[
\mathcal{R}_C(g^t) =  {\rm fix}_\T(h^t) = \bigcup_{w \in W} {\rm
fix}_\T(H,w).
\]
\end{teorema}
\begin{prova}
By the connectedness of ${\rm fix}_\T(H,w)$ we only need to prove
that each ${\rm fix}_\T(H,w)$ is chain recurrent. Let $g^t = e^t
u^t h^t$ be the Jordan decomposition of $g^t$. Note that the
restriction of $g^t$ to ${\rm fix}_\T(H,w)$ is given by $e^tu^t$.
First we show that for each $x \in \F_\T$ there exists $y \in
\F_\T$ such that $u^t x\to y$, when $t \to \pm\infty$. In fact,
$i(u^t x) = \rho(u)^t i(x)$.  By Lemma \ref{lemarepresent},
$\rho(u)$ is unipotent, so that, by Lemma \ref{lemaunipot}, there
exists $[v]$ such that $i(u^t x) \to [v]$, when $t \to \pm\infty$.
Using that $i$ is an embedding, there exists $y \in \F_\T$ such
that $i(y)=[v]$, which proves the claim. Now by Lemma
\ref{lemarepresent}, $\rho(e^t)$ is elliptic, so it lies in a
subgroup conjugated to $O(V)$.  This allows us to choose a metric
in $V$ such that $\rho(e^t)$ is an isometry for all $t \in \t$.
This metric induces a metric in ${\mathbb P} V$ and thus in
$\F_\T$, by using the Pl\"{u}cker embedding, so that $e^t$ is an
isometry in $\F_\T$. By Lemma \ref{lemaeu} applied to $u^t$, $e^t$
it follows that $g^t$ is chain recurrent on ${\rm fix}_\T(H,w)$.
\end{prova}

We remark that $H$ gives the parabolic type of $g^t$ and $Z_H$
gives the block reduction as defined in \cite{pss}.

Now we obtain the desired characterization of the recurrent set.

\begin{teorema}\label{teorecflag}
Let $g^t$ be a flow on $\F_\T$ and $g^t = e^th^tu^t$ its Jordan
decomposition. Then the recurrent set of $g^t$ in $\F_\T$ is given
by
\[
\mathcal{R}(g^t) = \emph{fix}_\T\left(h^t\right) \cap
\emph{fix}_\T\left(u^t\right).
\]
\end{teorema}
\begin{prova}
By Lemma \ref{lemarepresent} and by Theorem \ref{teorec}, we have
that
\[
\mathcal{R}(\rho(g)^t) = {\rm fix}(\rho(h)^t) \cap {\rm
fix}(\rho(u)^t).
\]
Thus the result follows by noting that
$$
i(\mathcal{R}(g^t)) = \mathcal{R}(\rho(g)^t) \cap i(\F_\T),
$$
$$
i( {\rm fix}_\T(h^t) ) = {\rm fix}(\rho(h)^t) \cap i(\F_\T)
\quad\mbox{and}\quad i( {\rm fix}_\T(u^t) ) = {\rm fix}(\rho(u)^t)
\cap i(\F_\T).
$$
\end{prova}

By using the previous characterization of the recurrent set, the
following result computes the topological entropy of linearly
induced flows on flag manifolds (see \cite{walters} for definition
and properties of topological entropy).

\begin{teorema}\label{teoentropy}
If $g^t$ is a flow in $\F_\T$, with $t \in \Z$, then its
topological entropy vanishes.
\end{teorema}
\begin{prova}
By using the variational principle and Poincar\'e recurrence
theorem (see \cite{walters}), the topological entropy of $g^t$
coincides with the topological entropy of its restriction to the
closure of its recurrent set. By Theorem \ref{teorecflag}, we have
that the recurrent set of $g^t$ is closed. Now let $g^t =
e^th^tu^t$ be the Jordan decomposition of $g^t$. Using again
Theorem \ref{teorecflag}, it follows that the restrictions of
$g^t$ and $e^t$ to $\mathcal{R}(g^t)$ coincide. Arguing exactly as
in the proof of Theorem \ref{teomorseflag}, we can provide a
metric in $\F_\T$ such that $e^t$ is an isometry in $\F_\T$, for
every $t \in \Z$. Thus the restriction of $e^t$ to
$\mathcal{R}(g^t)$ is also an isometry and therefore its
topological entropy vanishes.
\end{prova}

\subsection{Conley index and structural stability}

In this section, we first compute the Conley index of a linearly
induced flow $g^t$ on the flag manifold $\F_\T$.  Then we
characterize the linearly induced flows which are structurally
stable.

We say that $g^t$ is a conformal flow if $u^t = 1$, for all $t \in
\t$, where $g^t = e^t h^t u^t$ is its Jordan decomposition. For each
Iwasawa decomposition $G = KAN$ such that $H \in \cl \a^+$, we
define the conformal subgroup given by the direct product $C_H = K_H
A_{\T(H)}$.  We say that the flow $g^t$ has a conformal reduction if
there exists a conformal subgroup such that $g^t \in C_H$, for all
$t \in \t$ (this is the kind of linear flow considered in
\cite{hermann}).

\begin{proposicao}\label{propredconforme}
The flow $g^t$ is conformal if and only if it has a conformal
reduction.
\end{proposicao}
\begin{prova}
First we work with $\t = \R$.  Assume that $g^t$ has a conformal
reduction. Thus there exists a choice of an Iwasawa decomposition
$G = KAN$ with $H \in \cl \a^+$ such that $g^t=\exp(tX) \in C_H$,
for all $t \in \t$. Deriving at $t=0$ we get that $X$ belongs to
the Lie algebra of $C_H$ which is given by the direct sum $\k_H
\oplus \a_{\T(H)}$. Thus we can decompose $X = \widehat{E} +
\widehat{H}$ where $\widehat{E} \in \k_H$ commutes with
$\widehat{H} \in \a_{\T(H)}$. Since $\ad(\wh{E})$ is
$B_{\theta}$-anti-symmetric and $\ad(\wh{H})$ is
$B_{\theta}$-symmetric, it follows that $X = \widehat{E} +
\widehat{H}$ is the Jordan decomposition so that the nilpotent
part $N$ of $X$ vanishes and thus $u^t = 1$, for all $t \in \t$.
Conversely, assuming that $u^t = \exp(tN) = 1$, for all $t \in
\t$, we have that $N=0$, so that $X = E + H$ is its Jordan
decomposition. By Lemma \ref{lemajordan}, it follows that there
exists an Iwasawa decomposition $\g =\k \oplus \a \oplus \n^+$
with $H \in \cl \a^+$ such that $E \in \k_H$ and $H \in \cl \a^+$.
Hence $g^t=\exp(tX) \in C_H$, for all $t \in \t$.

Now we work with $\t=\Z$. Suppose that $g^t$ has a conformal
reduction. Hence there exists a choice of an Iwasawa decomposition
$G = KAN$ with $H \in \cl \a^+$ such that $g \in C_H$. Thus we can
decompose $g = \widehat{e}\widehat{h}$ where $\widehat{e} \in K_H$
commutes with $\widehat{h} \in A_{\T(H)}$. Since $\wh{e}$ is a
$B_{\theta}$-isometry and $\wh{h}$ is $B_{\theta}$-positive, this
is the multiplicative Jordan decomposition of $g$ so that $u^t =
1$, for all $t \in \t$. Conversely, assume that $u^t = 1$ so that
$g = eh$ is its multiplicative Jordan decomposition. By Lemma
\ref{lemajordan}, it follows that there exists an Iwasawa
decomposition $G = KAN$ with $H \in \cl \a^+$ such that $e \in
K_H$ and $H \in \cl \a^+$. Therefore $g^t=\exp(tX) \in C_H$, for
all $t \in \t$.
\end{prova}

Given an Iwasawa decomposition of $\g$ we recall that we can embed
the flag manifold $\F_\T$ in $\s$ in the following way.  Take
$H_{\Theta }\in \mathrm{cl}\frak{a}^{+}$ such that $\Theta =\Theta
(H_{\Theta })$ and put
$$
i: \F_\T \to \s, \qquad g \p_\T \mapsto k H_\T,
$$
where $g \in G$ and $g= kan$ is its Iwasawa decomposition, with $k
\in K$, $a \in A$, $n \in N$ (see Proposition 2.1 of \cite{dkv}).

\begin{proposicao}
Let the flow $g^t$ be conformal. Then there exists an Iwasawa
decomposition such that the height function of $H$ with respect to
the above embedding given by
$$
f: \F_\T \to \R, \qquad x \mapsto B_\theta(i(x),\, H),
$$
is a Lyapunov function for the finest Morse decomposition.
\end{proposicao}
\begin{prova}
Since $g^t$ is conformal, by the previous result there exists, an
Iwasawa decomposition such that $g^t \in C_H$. Decompose $g^t =
e^t h^t$ with $e^t \in K_H$ and $h^t \in A_{\T(H)}$.  Note that
height function $l$ is $K_H$-invariant since for $k \in K_H$, we
have
$$
f( kx ) = B_\theta(ki(x),\, H) = B_\theta(i(x),\, k^{-1}H) =
B_\theta(i(x),\, H) = f(x),
$$
where we used that $k$ is an isometry with respect to $B_\theta$.
Thus we have that
$$
f(g^t x) = f( h^t x),
$$
for all $t \in \t$, and the result follows by Proposition 3.3 item
(ii) of \cite{dkv}.
\end{prova}

The next result gives further information about the finest Morse
decomposition (see Theorem 3.2, Proposition 3.5 and Proposition
7.1 of \cite{pss}).  For the definition of the flag manifold
$\F(\Delta)_{H_0}$ see Section 3.3 of \cite{pss}.

\begin{proposicao}\label{proplinearizacao}
Let $g^t$ be a flow on $\F_\T$. Put $\Delta =\Theta (H)$ and take
$H_{\Theta }\in \mathrm{cl}\frak{a}^{+}$ such that $\Theta =\Theta
(H_{\Theta })$. Then ${\rm fix}_\T(H,w)$ is diffeomorphic to the
flag manifold $\mathbb{F}(\Delta)_{H_{0}}$, where $H_{0}$ is the
orthogonal projection of $wH_{\Theta }$ in $\frak{a}\left( \Delta
\right) $. Furthermore, the stable and unstable sets of ${\rm
fix}_\T(H,w)$ are diffeomorphic to vector bundles over
$\mathbb{F}(\Delta)_{H_{0}}$. Moreover, if $g^t$ is a conformal
flow, then it is normally hyperbolic and its restriction to
$\mathrm{st}_{\Theta }(H,w)$ and $\mathrm{un}_{\Theta }(H,w)$ are
conjugated to linear flows.
\end{proposicao}

In the next result we obtain the Conley index of the attractor and,
when $g^t$ is conformal, of all Morse components (see Proposition
7.6 and Corollary 7.8 of \cite{pss}).

\begin{teorema}
If $\T(H) \subset \T$ or $\T \subset \T(H)$, then the Conley index
of the attractor ${\rm fix}^+_\T(H)$ is the homotopy class of the
flag manifold $\F(\Delta)_{H_0}$.

Let $g^t$ be a conformal flow on $\F_\T$ and consider the same
notation of the previous proposition. Then the Conley index of the
Morse component ${\rm fix}_\T(H,w)$ is the homotopy class of the
Thom space of the vector bundle $\mathrm{un}_{\Theta }(H,w) \to
\mathbb{F}(\Delta)_{H_{0}}$. In particular, we have the following
isomorphism in cohomology
\[
CH^{*+n_{w}}({\rm fix}_\T(H,w))\simeq H^{*}\left(
\mathbb{F}(\Delta)_{H_{0}}\right),
\]
where $n_{w}$ is the dimension of $\mathrm{un}_{\Theta }(H,w)$ as a
vector bunble. The cohomology coefficients are taken in ${\mathbb
Z}_{2}$ in the general case and in ${\mathbb Z}$ when
$\mathrm{un}_{\Theta }(H,w)$ is orientable.
\end{teorema}

An element $H \in \a$ is said to be regular if there is no root
$\alpha \in \Sigma$ such that $\alpha(H)=0$. An element $X \in \g$
is h-regular if $H$ is regular in $\a$, for some Iwasawa
decomposition $\g = \k \oplus \a \oplus \n^+$, where $X = E + H +
N$ is its Jordan decomposition.  Note that if $X$ is h-regular
then so is $\psi X$ for $\psi \in \Int(\g)$. An element $g \in G$
is h-regular if $H$ is regular in $\a$, for some Iwasawa
decomposition $\g = \k \oplus \a \oplus \n^+$, where $g = ehu$ is
its Jordan decomposition and $\log h = \ad(H)$. Note that if $g$
is h-regular, then so is $hgh^{-1}$ for $h \in G$.

\begin{proposicao}\label{densidadereg}
The h-regular elements of $\g$ are dense in $\g$, and the
h-regular elements of $G$ are dense in $G$.
\end{proposicao}
\begin{prova}
Fix a Cartan involution $\theta$.  Let $\frak{j}$ be a
$\theta$-stable Cartan subalgebra. Take a Cartan decomposition
$\frak{j} = (\k \cap {\frak j}) \oplus (\s \cap {\frak j})$. Since
$\frak{j}$ is maximal abelian in $\g$, it follows that
$\frak{a}_{\frak{j}} = \s \cap {\frak j}$ is maximal abelian in
$\s$.  If $X \in \frak{j}$, then its Cartan and Jordan
decompositions coincide.  In fact, $\frak{j}$ is abelian, an
element in $\k$ is $B_{\theta}$-anti-symmetric and an element in
$\s$ is $B_{\theta}$-symmetric. Now since the regular elements of
$\frak{a}_{\frak{j}}$ are dense in $\frak{a}_{\frak{j}}$, it
follows that the h-regular elements of $\frak{j}$ are dense in
$\frak{j}$. By Proposition 1.3.4.1 p.101 of \cite{w} there exist
$\theta$-stable Cartan subalgebras $\frak{j}_1, \ldots,
\frak{j}_r$ such that the set
$$
\bigcup_{i=1}^r \Int(\g) \frak{j}_i
$$
is dense in $\g$.  Thus it follows that the h-regular elements of
$\frak{g}$ are dense in $\frak{g}$.

Let $J$ be a $\theta$-stable Cartan subgroup. By Proposition
1.4.1.2 p.109 of \cite{w}, we can take the Cartan decomposition
\[
J=(J\cap K)(\exp(\frak{j}\cap\frak{s})).
\]
Since $\frak{j}$ is maximal abelian in $\g$, it follows that
$\frak{a}_{\frak{j}} = \s \cap {\frak j}$ is maximal abelian in
$\s$.  If $g \in J$, then its Cartan and Jordan decompositions
coincide. In fact, $J$ centralizes $\exp(\s)$, an element in $K$
is a $B_{\theta}$-isometry and an element in $\exp(\s)$ is
$B_{\theta}$-positive. Now since the regular elements of
$\frak{a}_{\frak{j}}$ are dense in $\frak{a}_{\frak{j}}$, it
follows that the h-regular elements of $J$ are dense in $J$. By
Theorem 1.4.1.7 p.113 of \cite{w} there exist $\theta$-stable
Cartan subgroups $J_1, \ldots, J_r$ such that the set
$$
\bigcup_{i=1}^r \bigcup_{g \in G} g J_i g^{-1}
$$
is dense in $G$.  Thus it follows that the h-regular elements of
$G$ are dense in $G$.
\end{prova}

Note that a regular element is automatically conformal since we have
that $Z_H = MA$ when $H$ is regular, where $M$ is the centralizer of
$\a$ in $K$. It follows that the above result implies Theorem 8.1 of
\cite{hermann}.

\begin{lema}\label{lemanotreg}
Let $H \in \cl \a^+$ and $\T \subset \Sigma$.  If $H$ is not
regular and $\T \neq \Sigma$ then there exists $w \in W$ such that
${\rm fix}_\T(H,w)$ is not an isolated point in $\F_\T$.
\end{lema}
\begin{prova}
Given $w \in W$, by Proposition 3.6 p.326 of \cite{dkv} the map $
\varphi_w: w\n^-_\T \to N^-_\T b_\T \subset \F_\T$, $Y \mapsto
\exp(Y) wb_\T $ is a diffeomorphism such that $\varphi_w(0) =
wb_\T$ and
$$
\til{H} \circ \varphi_w = {\rm d} \varphi_w( \ad(H)|_{w\n^-_\T}),
$$
where $\til{H}(x)$ denotes the induced field of $H$ at $x$.  Since
${\rm fix}_\T(H,w)$ consists of the connected set of zeroes of the
induced vector field $\til{H}$ which pass trough $wb_\T$, it
follows that
$$
\varphi_w^{-1}( {\rm fix}(H,w)_\T \cap N^-_\T b_\T ) = {\rm Ker\,}
\ad(H)|_{w\n^-_\T}.
$$
Now if $H$ is not regular then there exists $\alpha \in \Sigma$
such that $\alpha(H) = 0$.  We have that
$$
w\n^-_\T = \sum \{ \g_{w^*\beta}:\, \beta \in \Pi^- - \langle \T
\rangle \}.
$$
Since $\T \neq \Sigma$ we can take $\beta \in \Pi^- - \langle \T
\rangle$.  Since the Weyl group $W$ acts transitively on $\Pi$ we
can take $w \in W$ such that $w^*\beta = \alpha$. Since $\alpha(H)
= 0$ it follows that $\g_\alpha \subset {\rm Ker\,}
\ad(H)|_{w\n^-_\T}.$  From the above discussion it follows that
$$
\varphi_w(\g_\alpha) \subset {\rm fix}_\T(H,w),
$$
so that ${\rm fix}_\T(H,w)$ is not an isolated point in $\F_\T$.
\end{prova}

Now we obtain the desired characterization of the linearly induced
flows which are structurally stable. For the concepts of
structural stability and of Morse-Smale flows and diffeomorphisms
see \cite{smale}. A flow $g^t$ in $\F_\T$ is regular if $H$ is
regular in $\a$, where $g^t = e^th^tu^t$ is its Jordan
decomposition.

\begin{teorema}\label{teoestab}
Let $g^t$ be a flow on $\F_\T$, where $\T \subset \Sigma$ with $\T
\neq \Sigma$. Then the following conditions are equivalent:
\begin{itemize}
\item[(i)] $g^t$ is regular,

\item[(ii)] $g^t$ is Morse-Smale and

\item[(iii)] $g^t$ is structurally stable.
\end{itemize}
\end{teorema}
\begin{prova}
First we show that condition (i) implies condition (ii). Let $g^t
= e^th^tu^t$ be the Jordan decomposition of $g^t$ and $H = \log
h$. If $H$ is regular, then for each $w \in W$, the set ${\rm
fix}_\T(H,w)$ reduces to a point. Thus, by Theorem
\ref{teomorseflag}, it follows that $\mathcal{R}_C(g^t)$ is a
finite set of fixed points. By Proposition \ref{proplinearizacao},
we have that each point in $\mathcal{R}_C(g^t)$ is hyperbolic.
Furthermore, by Lemma 4.2 of \cite{dkv} p.331, the stable and
unstable manifolds, given by Proposition \ref{propmorseflag},
intersect transversally. Thus it follows that $g^t$ is
Morse-Smale. By a result of \cite{smale}, it follows that
condition (ii) implies condition (iii).

Now we show that the negation of condition (i) implies the
negation of condition (iii). If $H$ is not regular, for each $\T
\subset \Sigma$ with $\T \neq \Sigma$, there exists a Morse
component ${\rm fix}_\T(H,w)$ which has infinite points. Thus, by
Theorem \ref{teomorseflag}, the same is true for
$\mathcal{R}_C(g^t)$. When $\t = \Z$, let $\widehat{g} \in G$ be a
h-regular element arbitrarily close to $g$, given by Proposition
\ref{densidadereg}, so that $\widehat{g}^t$ is a regular flow.
Then $\widehat{g}$ is arbitrarily close to $g$ in ${\rm
Diff}(\F_\T)$. By the first part of the proof, we have that
$\mathcal{R}_C(\widehat{g}^t)$ is finite and thus $\widehat{g}^t$
and $g^t$ cannot be topologically equivalent. This proves that, in
this case, $g^t$ is not structurally stable. In the case $\t =
\R$, we have that $g^t = \exp(tX)$, for some $X \in \g$. Let
$\widehat{X} \in \g$ be an h-regular element arbitrarily close to
$X$, given by Proposition \ref{densidadereg}, so that
$\widehat{g}^t = \exp(t\widehat{X})$ is a regular flow. Then
$\widehat{X}$ is arbitrarily close to $X$ in ${\mathcal
X}(\F_\T)$. Using again the first part of the proof, it follows
that $\mathcal{R}_C(\widehat{g}^t)$ is finite and thus
$\widehat{g}^t$ and $g^t$ cannot be topologically equivalent. Thus
we have again that $g^t$ is not structurally stable.
\end{prova}

\begin{obs}
Here we show how the results on the flag manifolds recover the
analogous ones on the projective space. Consider the Lie algebra
${\frak g} = {\rm sl}(n,\R)$ and the Lie group $G = \Int(\g) =
\Ad({\rm Sl}(n,\R))$.  A canonical choice of Iwasawa decomposition
gives the maximal abelian $\a = $ diagonal matrices with trace
zero and roots $\Pi$ given by the functionals in $\a$ defined by
$\alpha_{ij}({\rm diag}(H_k)) = H_i - H_j$, $i\neq j$, $i,j =
1,\ldots,n$.  Fixing the set of simple roots $\T =
\{\alpha_{i,i+1}:\, i=2,\ldots,n-1\}$, then the corresponding
parabolic subgroup in ${\rm Sl}(n,\R)$ is
\[
P_\T = \left( \textstyle{
\begin{array}{ccc}
a  & v \\
0  & A \\
\end{array}} \right),\, a\in\R,\, v\in\R^{n-1}\mbox{ and } {\rm
tr}A + a = 0,
\]
which is precisely the isotropy subgroup at $[e_1]$ of the
canonical action of ${\rm Sl}(n,\R)$ in ${\mathbb P}\R^n$.  It
follows that the map
\[
\varphi: {\mathbb P}\R^n \to \F_\T,\quad T[e_1] \mapsto \Ad(T)
\p_\T, \quad T \in {\rm Sl}(n,\R),
\]
is an equivariant diffeomorphism.  For $X \in \g$, the map
$\varphi$ conjugates the action of ${\rm e}^{tX}$ on ${\mathbb
P}\R^n$ with the action of $\exp(tX)$ on $\F_\T$, since
\[
\varphi( {\rm e}^{tX} T [e_1] ) = {\rm e}^{t\ad(X)} \varphi( T
[e_1] ) = \exp(tX) \varphi( T [e_1] ).
\]
We claim that $X \in {\rm sl}(n,\R)$ is h-regular if the real part
of its eigenvalues are distinct.  In fact, writing the Jordan
decomposition $X = E + H + N$, $H \in \cl \a^+$, then the entries
of $H = {\rm diag}(H_1,\ldots,H_n)$ are the real part of the
eigenvalues of $X$.  The claim follows, since $H$ is regular in
$\a$ when $\alpha_{ij}(H) = H_i - H_j \neq 0$.
\end{obs}

\section{Floquet theory}
In the previous sections, from the point of view of differential
equations, we have considered equations with constant
coefficients. In this section we extend these results to equations
with periodic coefficients.  Throughout this section we fix $\t =
\R$.

The fundamental solution associated with a given continuous map $t
\in \t \mapsto X(t) \in \g$ is the map $t \in \t \mapsto g(t) \in
G$ which satisfies
\[
g'(t) = X(t)g(t)
\]
and $g(0) = I$. It is straightforward to show that
\[
\rho_s : \t \to G, \qquad t \mapsto \rho_s(t) = g(t + s)g(s)^{-1}
\]
is the fundamental solution associated with the map $t \in \t
\mapsto X(s + t) \in \g$. When $X$ is constant, we have that $g(t)
= g^t = \exp(tX)$, which is the flow induced on $G$ by the
right-invariant vector field $X^r(a) = Xa$, where $a \in G$.

When $X$ is $T$-periodic $g(t)$ is not in general a flow. Despite
of this, we can associate a flow to $g(t)$ in the following way.
Since $X(t + T) = X(t)$, we have that $\rho_T(t) = g(t)$, for all
$t \in \t$. Thus it follows that $g(t + T) = g(t)g(T)$ and that
$g(t + mT) = g(t)g(T)^m$.  We will need the following result which
can be regarded as a refinement of Floquet's lemma to semisimple
Lie groups (see Theorem 2.47, p.163 of \cite{chicone}).

\begin{lema}
Let $G=\Int(\g)$. For $g \in G$ there exist $m \in \N$ and $X \in
\g$ such that $g^m = \exp(X)$.
\end{lema}
\begin{prova}
Let $g=ehu$ be its Jordan decomposition.  By Lemma \ref{lemajordan}
there exists an Iwasawa decomposition $G = KAN$ such that $h = {\rm
e}^{\ad(H)}$, $H \in \cl \a^+$ and $e \in K_H$. We have that $u = I
+ T$ where $T$ is a nilpotent map. Since $e$, $h$ commute with $u$,
it follows that $e$, $h$ commute with $T$. By Lemma IX.7.3 p.431 of
\cite{helgason}, we have that $u = {\rm e}^{\ad(N)}$ where $N \in
\g$ is such that $\ad(N)$ is nilpotent. By Lemma VI.4.5 p.270 of
\cite{helgason}, $\ad(N) = \log(u) = \log(I + T)$ which is a
polynomial in $T$.  It follows that $e$, $h$ commute with $\ad(N)$
and thus $e \in K_N$.  We claim that $H$ commutes with $N$.  Note
that if $p(x)$ is a polynomial, then $p(h)$ commutes with $\ad(N)$.
There exists a basis such that $\ad(H)$ and $h={\rm e}^{\ad(H)}$ are
given respectively by the diagonal matrices ${\rm
diag}(\lambda_1,\ldots,\lambda_n)$ and ${\rm diag}({\rm
e}^{\lambda_1},\ldots,{\rm e}^{\lambda_n})$. There exists $a,b > 0$
such that interval $[a,b]$ contains the eigenvalues of $h$.  By the
Weistrass approximation theorem, there exists a sequence of
polynomials $p_n(x)$ such that for $x \in [a,b]$ we have $p_n(x) \to
\log(x)$.  Thus we have that $p_n(h) \to \ad(H)$, which shows that
$\ad(H)$ commutes with $\ad(N)$, since $p_n(h)$ commutes with
$\ad(N)$ for all $n \in \N$.

From the above considerations it follows that $e$ lies in the
compact group $L = K_H \cap K_N$ which has Lie algebra $\k_H \cap
\k_N$.  It follows that there exists $m \in \N$ such that $e^m$ in
the connected component $L$ containing the identity.  Thus, by
Lemma II.6.10 p.135 of \cite{helgason}, there exist $E \in \k_H
\cap \k_N$ such that $e^m = \exp(E)$.  Taking $X = E + mH + mN$ it
follows that
\[
\exp(X) = \exp(E)\exp(H)^m\exp(N)^m = e^m h^m u^m = g^m.
\]
\end{prova}

By the above result, there exists $X \in \g$ such that $g(T)^m =
\exp(mTX)$. Defining $g^t = \exp(tX)$ and
\[
a(t) = g(t)g^{-t},
\]
it is straightforward to check that $a(t + mT) = a(t)$ and that
\[
\rho_s(t) = a(t + s)g^t a(s)^{-1},
\]
for all $t \in \t$. We have that the map defined by
\[
\phi\,^t(s, a) = (s + t, \rho_s(t)a)
\]
is a flow of automorphisms on the principal bundle $S^1 \times G$.
In fact, observing that the map $(s, a) \mapsto (s, a(s)a)$ is a
diffeomorphism of $S^1 \times G$ onto itself, it follows that
\[
\phi\,^t(s, a(s)a) = (s + t, a(s + t)g^ta)
\]
and thus is not hard to check that $\phi\,^t$ is already a flow.

From now on we consider the Jordan decomposition $g^t = e^th^tu^t$
of the flow on $\F_\T$ associated to $\phi\,^t$, where $\log h =
\ad(H)$ and $H$ lies in the closure of a fixed Weyl chamber $\a^+$
(see Section \ref{subsecjordan}). For each $\F_\T$, we can induce
a flow on $S^1 \times \F_\T$ by simply putting
\[
\phi\,^t(s, x) = (s + t, \rho_s(t)x)
\]
or
\[
\phi\,^t(s, a(s)x) = (s + t, a(s + t)g^tx).
\]
Note that when $X(t)$ is constant and equals to $X$, then $a(t) =
I$ and $\phi\,^t(s, x) = (s + t, g^tx)$, so we return to the same
context of the previous section.

Now we obtain the desired characterization of the recurrent set.

\begin{teorema}
The recurrent set of $\phi\,^t$ in $S^1 \times \F_\T$ is given by
\[
\mathcal{R}(\phi\,^t) = \{(s, a(s)x) :\, s \in S^1,\,  x \in
\emph{fix}_\T\left(h^t\right) \cap
\emph{fix}_\T\left(u^t\right)\}.
\]
\end{teorema}
\begin{prova}
Denoting by
$$
R = \{(s, a(s)x) :\, s \in S^1,\,  x \in {\rm
fix}_\T\left(h^t\right) \cap {\rm fix}_\T\left(u^t\right)\},
$$
we will first show that $R \subset \mathcal{R}(\phi\,^t)$.  Given
$(s,a(s)x)$, where $x \in \mathcal{R}(g^t)$, by Theorem
\ref{teorecflag}, we have that $g^tx = e^tx$, for all $t \in \t$.
Arguing exactly as in the proof of Theorem \ref{teomorseflag}, we
can provide a metric $d$ in $\F_\T$ such that $e^t$ is an isometry
in $\F_\T$, for every $t \in \t$. By the compactness of $\F_\T$
and by Lemma \ref{lemaiso}, there exists a sequence $n_k\to
\infty$ such that $g^{n_kmT}x \to x$.  It follows that
\[
\phi\,^{n_kmT}(s, a(s)x) = (n_kmT + s, a(n_kmT + s)g^{n_kmT}x) =
(s, a(s)g^{n_kmT}x) \to (s,a(s)x),
\]
showing that $R \subset \mathcal{R}(\phi\,^t)$.  Conversely, let
$(s,a(s)x) \in \mathcal{R}(\phi\,^t)$.  Thus there exists $t_k \to
\infty$ such that
\[
(t_k + s, a(t_k + s)g^{t_k}x) = \phi\,^{t_k}(s, a(s)x)  \to
(s,a(s)x).
\]
Therefore $t_k + s \to s$ modulo $mT$ so that $a(t_k + s) \to
a(s)$ and thus $g^{t_k}x \to x$.
\end{prova}

Let $\{{\rm fix}_\T(H,w): w \in W_H\backslash W/W_\T \}$ be the
finest Morse decomposition of the flow $g^t$ given by Theorem
\ref{teomorseflag} and define
\[
{\cal M}_\T(H,w) = \{(s, a(s)x) :\, s \in S^1,\,  x \in {\rm
fix}_\T(H,w)\},
\]
which is a $\phi\,^t$-invariant subset of $S^1 \times \F_\T$. If $f
: \F_\T \to \R$ is a Lyapunov function for the finest Morse
decomposition of $g^t$ and defining
\[
F(s, a(s)x) = f(x),
\]
we have that $F : S^1 \times \F_\T \to \R$ is a Lyapunov function
for the family
\[
\{{\cal M}_\T(H,w) : w \in W_H\backslash W/W_\T \},
\]
which is, therefore, a Morse decomposition of the flow $\phi\,^t$.
In fact,
\[
F(\phi\,^t(s, a(s)x)) = F(s + t, a(s + t)g^tx) = f(g^tx)
\]
and thus $F \circ \phi\,^t$ is constant over each ${\cal
M}_\T(H,w)$ and strictly decreasing out of their union. Now we
characterize the stable and unstable sets of these Morse
components.

\begin{proposicao}\label{propmorsefloquet}
The stable and unstable sets of ${\cal M}_\T(H,w)$ are given by
$$
{\rm st} ({\cal M}_\T(H,w)) = \{(s, a(s)x) :\, s \in S^1,\,  x \in
{\rm st}_\T(H,w) \}
$$
and
$$
{\rm un} ({\cal M}_\T(H,w)) = \{(s, a(s)x) :\, s \in S^1,\,  x \in
{\rm un}_\T(H,w) \}.
$$
\end{proposicao}
\begin{prova}
Taking $x \in {\rm st}_\T(H,w)$, by Proposition
\ref{propmorseflag}, we have that $g^tx \to {\rm fix}_\T(H,w)$.
Thus it follows that
\[
\phi\,^t(s, a(s)x) = (s + t, a(s + t)g^tx) \to {\rm st} ({\cal
M}_\T(H,w)),
\]
showing that $(s, a(s)x) \in {\rm st} ({\cal M}_\T(H,w))$. The
equality follows by observing that, by the Bruhat decomposition,
the sets ${\rm st}_\T(H,w)$, $w \in W$, exhaust $\F_\T$.  The
proof for the unstable set is entirely analogous.
\end{prova}

We denote the only attractor and the only repeller ${\cal
M}_\T(H,1)$ and ${\cal M}_\T(H,w^-)$, respectively, by ${\cal
M}^+_\T(H)$ and ${\cal M}^-_\T(H)$. Using the previous result, we
obtain the desired characterization of the finest Morse
decomposition.

\begin{teorema}\label{teomorsefloquet}
 Each ${\cal M}_\T(H,w)$ is chain
transitive, so that $\{{\cal M}_\T(H,w) : w \in W_H\backslash
W/W_\T \}\}$ is the finest Morse decomposition of $\phi\,^t$. In
particular, the chain recurrent set of $\phi\,^t$ in $S^1 \times
\F_\T$ is given by
\[
\mathcal{R}_C(\phi\,^t) =  \{(s, a(s)x) :\, s \in S^1,\,  x \in
\emph{fix}_\T\left(h^t\right)\}.
\]
\end{teorema}
\begin{prova}
First assume that, for all given $\varepsilon >0$, $t \in \t$, and
$x, y \in {\rm fix}_\T(H,w)$, there is an $(\varepsilon,t)$-chain
from $(0,x)$ to $(s,a(s)y)$. With this we will construct an
$(\varepsilon,t)$-chain from $(\wh{s},a(\wh{s})x)$ to
$(\wt{s},a(\wt{s})y)$, for all $x, y \in {\rm fix}_\T(H,w)$.  In
fact, take $\delta > 0$ given by the $\varepsilon$-uniform
continuity of $\phi\,^{\wh{s}}$ in $S^1\times \F_\T$.  Denote by
$\ov{d}$ the metric in $S^1\times \F_\T$ given by
$$
\ov{d}((s,x),(r,y)) = |s-r| + d(x,y),
$$
where $d$ is a metric in $\F_\T$. Consider the $(\delta,t)$-chain
from $(0,x)$ to
$$
\left(\wt{s}-\wh{s},a(\wt{s}-\wh{s})g^{-\wh{s}}y\right),
$$
given by $t_i > t$ and $\eta_i \in S^1 \times \F_\T$, where
$i=1,\ldots,n+1$ such that $\ov{d}( \eta_i,\, \phi\,^{t_i}(\eta_i)
) < \delta,$ for $i=1,\ldots,n$. We have thus that, with the same
$t_i > t$, $\phi\,^{\wh{s}}(\eta_i)$ is an $(\varepsilon,t)$-chain
from $(\wh{s},a(\wh{s})x)$ to $(\wt{s},a(\wt{s})y)$.

Now we will prove the above assumption.  Let $x, y \in {\rm
fix}_\T(H,w)$ and take $\delta > 0$ given by the
$\varepsilon$-uniform continuity of the map $(s,z) \mapsto a(s)z$
in $S^1 \times \F_\T$. By the compactness of $\F_\T$, there exists
$\tau > 0$ such that $d(g^t z,z) < \delta/2$, for all $t \in
[0,\tau]$ and all $z \in \F_\T$. By Theorem \ref{teomorseflag},
there exists a $(\delta/2,t)$-chain from $x$ to $y$ in $\F_\T$
given by $t_i > t$ and $x_i \in \F_\T$, where $i=1,\ldots,n+1$.
Note that $n$ can be taken arbitrarily large such that $mT/n <
\tau$. Let $l$ be such that
$$
\wh{s} = s + lmT - (t_1 + \cdots + t_n) \in [0,mT].
$$
Consider $\wh{t} = \wh{s}/n \in [0,\tau]$ and define $\wh{t}_i =
t_i + \wh{t} > t$,
$$
\xi_1 = (0,x) \qquad \mbox{and} \qquad \xi_{i+1} = \left(\wh{t}_1
+ \cdots + \wh{t}_i,\, a(\wh{t}_1 + \cdots +
\wh{t}_i)x_{i+1}\right).
$$
We claim that this provides an $(\varepsilon,t)$-chain from
$(0,x)$ to $(s,a(s)y)$. In fact, note that
$$
\wh{t}_1 + \cdots + \wh{t}_n = s + lmT,
$$
and that $x_{n+1} = y $, so that we have $\xi_{n+1} = (s,a(s)y)$.
Since
$$
\phi\,^{\wh{t}_i}(\xi_i) = \left(\wh{t}_1 + \cdots + \wh{t}_i,\,
a(\wh{t}_1 + \cdots + \wh{t}_i) g^{\wh{t}_i} x_{i}\right)
$$
and
$$
d\left(g^{\wh{t}_i} x_{i},x_{i+1}\right) \leq d\left(g^{\wh{t}}z,
z\right) + d\left(g^{t_i}x_i, x_{i+1}\right) < \delta,
$$
where $z = g^{t_i} x_{i}$, we have that
$$
\ov{d}\left(\phi\,^{\wh{t}_i}(\xi_i), \, \xi_{i+1}\right) <
\varepsilon.
$$
\end{prova}

Consider the set
$$
Q_\phi = \{(s, a(s)a) :\, s \in S^1,\, a \in Z_H \} \subset S^1
\times G.
$$
Since $g^t \in Z_H$, it follows that $Q_\phi$ is a
$\phi\,^t$-invariant $Z_H$-reduction of the principal bundle $S^1
\times G$ which in \cite{pss} has been called a block reduction of
the flow $\phi\,^t$.  Looking at $S^1 \times \F_\T$ as an
associated bundle of $S^1\times G$, by the previous results, this
reduction has the following immediate properties
\[
{\cal M}_\T(H,w) =  Q_\phi \cdot w b_\T \qquad\mbox{and}\qquad
{\rm st} ({\cal M}_\T(H,w)) = Q_\phi \cdot {\rm st}_\T(H,w),
\]
recovering, in this context, Theorem 4.1 item (i) and Theorem 5.3
of \cite{pss}.

We say that the flow $\phi\,^t$ is conformal if its associated
flow $g^t$ is conformal.  In this case, by Proposition
\ref{propredconforme}, there exists a conformal subgroup $C_H$ of
$Z_H$ which contains $g^t$. It follows that the set
$$
C_\phi = \{(s, a(s)c) :\, s \in S^1,\, c \in C_H \} \subset S^1
\times G
$$
is a $\phi\,^t$-invariant $C_H$-reduction of the principal bundle
$S^1 \times G$ which in \cite{pss} has been called a conformal
reduction of the flow $\phi\,^t$.  Note that both $Q_\phi$ and
$C_\phi$ have a global section given by $s \mapsto (s,a(s))$.
Thus, by Corollary 8.4, p.36 of \cite{steenrod}, these are trivial
bundles, implying that their associated bundles are also trivial.
The next result gives further information about the finest Morse
decomposition (see Propositions 4.2, 6.2 and 7.1 of \cite{pss}).
For the definition of the flag manifold $\F(\Delta)_{H_0}$ see
Section 3.3 of \cite{pss}.

\begin{proposicao}\label{proplinearizacaofloquet}
Put $\Delta =\Theta (H)$ and take $H_{\Theta }\in
\mathrm{cl}\frak{a}^{+}$ such that $\Theta =\Theta (H_{\Theta })$.
Then ${\cal M}_\T(H,w)$ is homeomorphic to $S^1 \times
\mathbb{F}(\Delta)_{H_{0}}$, where $H_{0}$ is the orthogonal
projection of $wH_{\Theta }$ in $\frak{a}\left( \Delta \right) $.
Furthermore, the stable and unstable sets of ${\rm fix}_\T(H,w)$
are diffeomorphic to vector bundles over ${\cal M}_\T(H,w)$.
Moreover, if $\phi\,^t$ is a conformal flow, then it is normally
hyperbolic and its restriction to ${\rm st} ({\cal M}_\T(H,w))$
and ${\rm un} ({\cal M}_\T(H,w))$ are conjugated to linear flows.
\end{proposicao}

In the next result we obtain the Conley index of the attractor and,
when $\phi\,^t$ is conformal, of all Morse components (see Theorem
7.4 and Corollary 7.8 of \cite{pss}).

\begin{teorema}
If $\T(H) \subset \T$ or $\T \subset \T(H)$, then the Conley index
of the attractor ${\cal M}^+_\T(H)$ is its homotopy class.

If $\phi\,^t$ is conformal, then the Conley index of the Morse
component ${\cal M}_\T(H,w)$ is the homotopy class of the Thom space
of the vector bundle $\mathrm{un}({\cal M}_\T(H,w)) \to {\cal
M}_\T(H,w)$. In particular, we have the following isomorphism in
cohomology
\[
CH^{*+n_{w}}({\cal M}_\T(H,w))\simeq H^{*}\left( S^1 \times
\mathbb{F}(\Delta)_{H_{0}} \right),
\]
where $n_{w}$ is the dimension of $\mathrm{un}({\cal M}_\T(H,w))$ as
a vector bunble. The cohomology coefficients are taken in ${\mathbb
Z}_{2}$ in the general case and in ${\mathbb Z}$ when
$\mathrm{un}_{\Theta }(H,w)$ is orientable.
\end{teorema}


It follows that the (co)homology Conley indexes can be computed by
K\"uneth formula, once we know the (co)homology of the real flag
manifolds.  For the homology of real flag manifolds see
\cite{kocherlakota}.

\appendix

\section{Dynamics in projective spaces}

In this appendix, we relate the Jordan decomposition of $g^t$ in
${\rm Gl}(V)$ to the dynamics of the induced linear flow $g^t$ on
the projective space ${\mathbb P}V$, where $V$ is a finite
dimensional vector space. The main results of the section deals with
the characterization of the recurrent set and the finest Morse
decomposition in terms of the fixed points of the Jordan components.
We will need some preliminary lemmas.

\begin{lema}\label{lemathiago1}
Let $V = U \oplus W$ and let $x_n = u_n + w_n$ be a sequence with
$u_n \in U$ and $w_n \in W$.  Suppose that $u_n \neq 0$ for all $n
\in \N$ and that $\lim \frac{w_n}{|u_n|} = 0$.  Passing to the
projective space, if $[x_n] \in {\mathbb P}V$ converges to $[x] \in
{\mathbb P}V$ then $[x] \in {\mathbb P}U$.
\end{lema}
\begin{prova}
Without loss of generality we can suppose that $|u_n| = 1$ and that
$\lim w_n = 0$ since $[x_n] = [\frac{u_n}{|u_n|} +
\frac{w_n}{|u_n|}]$ and $\lim \frac{w_n}{|u_n|} = 0$. Now since
$|u_n| = 1$ it has a convergent subsequence $u_{n_k} \to u$, where
$u \in U$, since subspaces are closed in $V$.  Then
$$
[x] = \lim_{k \to \infty} [x_{n_k}] = \lim_{k \to \infty} [u_{n_k} +
w_{n_k}] = [u]
$$
and thus $[x] \in {\mathbb P}U$.
\end{prova}

\begin{lema}\label{lemathiago2}
Fix a norm $|\cdot|$ in $V$.  If $h = I$ then for each $x \neq 0$
there exists $\epsilon > 0$ such that $|g^t x| > \epsilon$ for all
$t \in \t$.
\end{lema}
\begin{prova}
By the Jordan canonical form, in an appropriate basis, $u$ is upper
triangular with ones on the diagonal.  Write $x$ in this basis as $x
= (x_1, \ldots, x_k, 0,\ldots 0)$ where $x_k$ is the last nonzero
coordinate of $x$.  Then $u^t$ fixes the last coordinate $x_k$ of
$x$ so that, if we take the euclidian norm $|\cdot|_1$ relative to
this basis, we have that $|u^t x|_1 \geq |x_k|$ for all $t \in \t$.
Taking the norm $|\cdot|_2$ which makes $e$ an isometry, we have
that $|v|_2 \geq C |v|_1$ for all $v \in V$, where $C > 0$. Since by
hypothesis $g^t = e^t u^t$, we have that
$$
|g^t x|_2 = |u^t x|_2 \geq C |u^t x|_1 \geq C|x_k|
$$
for all $t \in \t$.  Using that $|\cdot|_2$ is equivalent to the
norm $|\cdot|$ fixed in $V$, the lemma follows.
\end{prova}

For the following lemma, we need to recall the definition of the
spectral radius of a linear map $g$, which is given by
$$
r(g) = \max \{|\lambda|: \lambda \mbox{ is an eigenvalue of }g \}.
$$
Let $V$ have a norm $|\cdot|$.  We also denote by $|g|$ the
corresponding operator norm.

\begin{lema}\label{lemalucas1}
If $r(g) < 1$ then $|g^t| \to 0$.
\end{lema}
\begin{prova}
Let $g^t = e^t h^t u^t$ be the Jordan decomposition of the linear
flow $g^t$. Since the norms in $V$ are equivalent, we can choose
$|\cdot|$ such that the eigenvector basis of $h$ is orthonormal.
Then we have that $|h^t| = r(g)^t$.  Since $e$ is elliptic, it lies
in a subgroup conjugated to $O(V)$ so that the norm $|e^t|$ is
bounded, say by $M > 0$.  Since $u$ is unipotent, we have that $u =
e^N$, where $N$ is nilpotent, so that
$$
u^t = e^{tN} = I + tN + \cdots + (tN)^k/k!,
$$
for some fixed $k \in \N$.  Using that $|N^l| \leq |N|^l$, it
follows that for $v \in V$
$$
|u^t v| \leq |v|(1 + |N|t + \cdots + (|N|^k/k!) t^k) = |v|p(t),
$$
where $p(t)$ is a polynomial in $t$, so that $|u^t| \leq p(t)$.
Collecting the above results, we have that
$$
|g^t| \leq |e^t||h^t||u^t| \leq M r(g)^t p(t) \to 0,
$$
when $t \to \infty$, since $r(g) < 1$.
\end{prova}

\begin{lema}\label{lemthiago}
Let $V = V_{\lambda_1} \oplus V_{\lambda_2} \oplus \cdots \oplus
V_{\lambda_n}$ be the eigenspace decomposition of $V$ associated to
$h$ where $\lambda_1 > \lambda_2 > \cdots \lambda_n
> 0$.  Let $v = v_1 + v_2 + \cdots + v_n$, $v \neq 0$, with $v_i \in
V_{\lambda_i}$. Take $i$ the first index $k$ with $v_k \neq 0$ and
$j$ the last index $k$ with $v_k \neq 0$.  Then $\omega([v]) \subset
{\mathbb P} V_i$ and $\omega^*([v]) \subset {\mathbb P} V_j$.
\end{lema}
\begin{prova}
Denote by $g_k$ the restriction of $g/\lambda_i$ to $V_{\lambda_k}$.
We have that $g_k$ has spectral radius $\lambda_k/\lambda_i$ and
that $g_i$ has hyperbolic part equal to the identity. By Lemma
\ref{lemalucas1} we have that for $k > i$ we have $|g_k^t v_k| \leq
|g_k^t||v_k| \to 0$, when $t \to \infty$. By Lemma \ref{lemathiago2}
we have that $|g_i v_i| \geq \epsilon$, for some $\epsilon > 0$. Now
let $t_j \to \infty$ be such that $\lim_{j \to \infty} g^{t_j}[v] =
[x]$, then
\[
g^{t_j}[v] = \left[ \frac{g^{t_j}}{\lambda_i}(v_i + \cdots + v_t)
\right] = \left[g_i^{t_j}v_i + \sum_{k>i} g_k^{t_j}v_k \right] \to
[x],
\]
so that, by the previous arguments and Lemma \ref{lemathiago1}, it
follows that $[x] \in {\mathbb P}V_i$.
\end{prova}

With the previous result we obtain that the projective eigenspaces
of the hyperbolic part of $g^t$ give are Morse components for the
flow $g^t$ in the projective space.

\begin{proposicao}\label{propmorse}
The set $\{{\mathbb P} V_{\lambda_1}, \ldots, {\mathbb P}
V_{\lambda_n}\}$ is a Morse decomposition.  Furthermore, the stable
set of ${\mathbb P} V_{\lambda_i}$ is given by
$$
{\rm st} ({\mathbb P}V_{\lambda_i}) = \{ [v_i + v_{i+1} + \cdots +
v_n]:\, v_i \neq 0\},
$$
where $v_k \in V_{\lambda_k}$.
\end{proposicao}
\begin{prova}
Since $h^t$ commutes with $g^t$ and taking $v_k \in V_{\lambda_k}$,
it follows that
$$
h^t g^t v_k = g^t h^t v_k = \lambda_k g^t v_k,
$$
showing that $ V_{\lambda_k}$ is $g^t$-invariant.  The proposition
then follows from the definition of Morse decomposition, using the
previous lemma.
\end{prova}

\begin{figure}
    \begin{center}
        \includegraphics[scale=0.75]{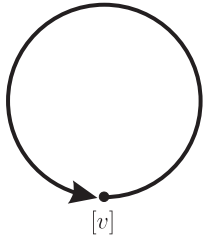}
    \end{center}
    \caption{\label{figura1}
    Unipotent element acting on the projective line.}
\end{figure}

In order to show that the above Morse decomposition is the finest
one, we need to consider the behavior of the unipotent component of
$g^t$.  This is done in the next lemma, which generalizes the
behavior on the projective line (see Figure \ref{figura1}) of the
action of the linearly induced map ${\rm e}^{tN}$, where
\[
N = \left( \textstyle{
\begin{array}{cc}
0 & 1 \\
0 & 0 \\
\end{array}} \right).
\]

\begin{lema}\label{lemaunipot}
Let $x \neq 0$ be a vector and $N$ be a nilpotent transformation. If
$k$ is such that $N^{k+1}x=0$ and $v = N^k x \neq 0$, then
$e^{tN}[x] \to [v]$, when $t \to \pm \infty$, where $t \in \t$.
Furthermore $e^{tN}v = v$, for all $t \in \t$.
\end{lema}
\begin{prova}
First note that
$$
e^{tN}v = \left(\sum_{j\geq 0} \frac{t^j}{j!}N^j\right)N^k x = v +
\sum_{j\geq 1} \frac{t^j}{j!}N^{k+j}x = v.
$$
Now we have that
$$
e^{tN}[x] = \left[x + tN + \cdots + \frac{t^k}{k!} v\right]
$$
and, multiplying by $k!/t^k$, we get that
$$
e^{tN}[x] = \left[v + \frac{k!}{t^k}\left(tNx + \cdots +
\frac{t^{k-1}}{(k-1)!} N^{k-1}x \right) \right] \to [v],
$$
when $t \to \pm \infty$, where $t \in \t$.
\end{prova}

Collecting the previous results we obtain the desired
characterization of the finest Morse decomposition.

\begin{teorema}\label{teoreccadeias}
Let $g : V \to V$ be a linear isomorphism, where $V$ is a finite
dimensional vector space.   Let $V = V_{\lambda_1} \oplus
V_{\lambda_2} \oplus \cdots \oplus V_{\lambda_n}$ be the eigenspace
decomposition of $V$ associated to $h$.  Then each ${\mathbb P}
V_{\lambda_i}$ is chain transitive, so that $\{{\mathbb P}
V_{\lambda_1}, \ldots, {\mathbb P} V_{\lambda_n}\}$ is the finest
Morse decomposition.  In particular, the chain recurrent set of $g$
in ${\mathbb P}V$ is given by
\[
\mathcal{R}_C(g^t) = \emph{fix}\left(h^t\right) = \bigcup_i {\mathbb
P} V_{\lambda_i}.
\]
\end{teorema}
\begin{prova}
By the connectedness of ${\mathbb P} V_{\lambda_i}$ we just need to
prove that each ${\mathbb P} V_{\lambda_i}$ is chain recurrent. We
note that the second equality on the equation above is immediate.
Thus, by Proposition \ref{propmorse} we have that $\mathcal{R}_C(g)
\subset {\rm fix}\left(h^t\right)$. Now we prove that ${\rm
fix}\left(h^t\right)$ is chain recurrent.  First note that the
restriction of $g^t$ to ${\rm fix}\left(h^t\right)$ is given by $e^t
u^t$. Since $e^t$ is elliptic, it lies in a subgroup conjugated to
$O(V)$.  This allows us to choose a metric in $V$ such that $e^t$ is
an isometry for all $t \in \t$. This metric induces a metric in
${\mathbb P} V$ such that $e^t$ is an isometry in ${\mathbb P} V$.
By Lemmas \ref{lemaunipot} and \ref{lemaeu} applied to $u^t$, $e^t$
it follows that $g^t$ is chain recurrent on ${\rm
fix}\left(h^t\right)$.
\end{prova}

We conclude with the desired characterization of the recurrent set.

\begin{teorema}\label{teorec}
Let $g : V \to V$ be a linear isomorphism, where $V$ is a finite
dimensional vector space. Then the recurrent set of $g$ in ${\mathbb
P}V$ is given by
\[
\mathcal{R}(g^t) = \emph{fix}\left(h^t\right) \cap
\emph{fix}\left(u^t\right).
\]
\end{teorema}
\begin{prova}
Let $[x]$ be such that $[x] \in \omega([x])$.  By Theorem
\ref{teoreccadeias} we have that $[x] \in {\rm
fix}\left(h^t\right)$.  Let $t_j \to \infty$ be such that
$g^{t_j}[x] \to [x]$.  Since $e^t$ is elliptic, it lies in a
subgroup conjugated to $O(V)$, so we can assume that $e^{t_j} \to
E$.  Note that $E$ commutes with the Jordan components of $g$. By
Lemma \ref{lemaunipot}, there exists a fixed point $[v]$ of $u^t$
such that $u^{t_j}[x] \to [v]$.  Since $g^{t_j} =
e^{t_j}u^{t_j}h^{t_j}$ it follows that
$$
[x] = \lim g^{t_j} [x] = \lim e^{t_j} u^{t_j} [x] = E [v].
$$
The theorem follows since $E$ commutes with $u^t$ and $[v]$ is a
$u^t$-fixed point.
\end{prova}

We illustrate the above results with some examples in dimension
three. In order to stay in the context of semisimple Lie groups we
work in ${\rm Sl}(3,\R)$.

\begin{exemplo}
Let $X \in {\rm sl}(3,\R)$.  There exists $g \in {\rm Sl}(3,\R)$
such that $gXg^{-1}$ has one of the following Jordan canonical
forms, where $a,b \in \R$:
\[
X_1 = \left( \textstyle{
\begin{array}{ccc}
-a & 0 & 0 \\
0 & -b & 0 \\
0 & 0 & a+b\\
\end{array}} \right),\quad
X_2 = \left( \textstyle{
\begin{array}{ccc}
0 & 1 & 0 \\
0 & 0 & 1 \\
0 & 0 & 0 \\
\end{array}} \right),\quad
X_3 = \left( \textstyle{
\begin{array}{ccc}
0 & 0 & 1 \\
0 & 0 & 0 \\
0 & 0 & 0\\
\end{array}} \right),
\]
\[
X_4 = \left( \textstyle{
\begin{array}{ccc}
-a & -b & 0 \\
b  & -a & 0 \\
0  & 0  & 2a\\
\end{array}} \right)\quad \mbox{and} \quad
X_5 = \left( \textstyle{
\begin{array}{ccc}
-a & 1 & 0 \\
0  & -a & 0 \\
0  & 0  & 2a\\
\end{array}} \right).
\]
Let $a, b > 0$. We have that the nilpotent component of $X_4$ is
zero, while its elliptic and hyperbolic components are given,
respectively, by
\[
E = \left( \textstyle{
\begin{array}{ccc}
0  & -b & 0 \\
b  & 0 & 0 \\
0  & 0  & 0\\
\end{array}} \right)\quad \mbox{and} \quad
H = \left( \textstyle{
\begin{array}{ccc}
-a & 0  & 0 \\
0  & -a & 0 \\
0  & 0  & 2a\\
\end{array}} \right).
\]
In this case, the recurrent and chain recurrent sets coincide and we
have two Morse components: the attractor $[e_3]$ and the repeller
$[\R e_1 \oplus \R e_2]$ (see Figure \ref{fig2}(a)). We also have
that the elliptic component of $X_5$ is zero, while its hyperbolic
and nilpotent components are given, respectively, by
\[
H = \left( \textstyle{
\begin{array}{ccc}
-a & 0  & 0 \\
0  & -a & 0 \\
0  & 0  & 2a\\
\end{array}} \right)\quad \mbox{and} \quad
N = \left( \textstyle{
\begin{array}{ccc}
0 & 1 & 0 \\
0 & 0 & 0 \\
0 & 0 & 0 \\
\end{array}} \right).
\]
In this case, the recurrent and chain recurrent sets are different,
but the Morse components remain the same.  The recurrent set is
given by $\{[e_1],\, [e_3]\}$ (see Figure \ref{fig2}(b)). In both
cases, the stable set of the attractor is the complement of the
repeller.

\begin{figure}
    \begin{center}
        \includegraphics[scale=0.75]{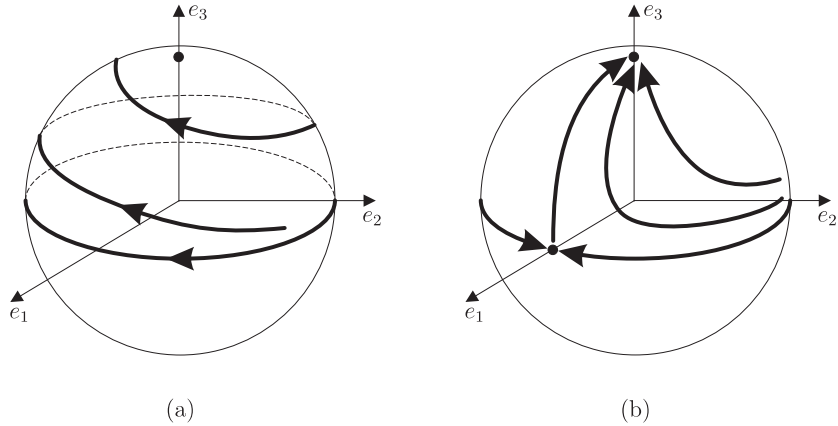}
    \end{center}
    \caption{\label{fig2}
    Dynamics on the projective space represented on the two-sphere.}
\end{figure}

\end{exemplo}

\end{document}